\documentclass[11pt,sfbold,
]{paper}
\usepackage{amsfonts}
\usepackage{amssymb}
\usepackage{amsmath}
\usepackage[USenglish]{babel}
\usepackage{theorem}
\usepackage{fancyhdr}
\usepackage[letterpaper]{geometry}
\usepackage{setspace}

\newtheorem{theorem}{Theorem}
\newtheorem{example}[theorem]{Example}

\newtheorem{corollary}[theorem]{Corollary}

\newtheorem{definition}[theorem]{Definition}

\newtheorem{lemma}[theorem]{Lemma}

\newtheorem{remark}[theorem]{Remark}

\input{tcilatex}

\begin{document}

\author{Edgar Delgado-Eckert\thanks{%
Centre for Mathematical Sciences, Munich University of Technology,
Boltzmannstr.3, 85747 Garching, Germany.} \thanks{%
Pathology Department, Tufts University, 150 Harrison Ave., Boston, MA 02111,
USA.} \thanks{%
Present affiliation (and correspondence address): ETH Z\"{u}rich, Department
of Biosystems Science and Engineering (D-BSSE), WRO-1058-7.42, Mattenstrasse
26, 4058 Basel, Switzerland. Email: edgar.delgado-eckert@mytum.de}}
\title{An algebraic and graph theoretical framework to study monomial
dynamical systems over a finite field}
\maketitle

\begin{abstract}
A monomial dynamical system $f:K^{n}\rightarrow K^{n}$ over a finite field $%
K $ is a nonlinear deterministic time discrete dynamical system with the
property that each component function $f_{i}:K^{n}\rightarrow K$ is a monic
nonzero monomial function. In this paper we provide an algebraic and graph
theoretic framework to study the dynamic properties of monomial dynamical
systems over a finite field. Within this framework, characterization
theorems for fixed point systems (systems in which all trajectories end in
steady states) are proved. In particular, we present an algorithm of
polynomial complexity to test whether a given monomial dynamical system over
a finite field is a fixed point system. Furthermore, theorems that
complement previous work are presented and alternative proofs to previous
results are supplied.
\end{abstract}

\section{Introduction}

Time discrete dynamical systems over a finite set $X$ are an important
subject of active mathematical research. One relevant example of such
systems are \textit{cellular automata}, first introduced in the late 1940s
by John von Neumann (e.g., \cite{MR0299409}). More general examples of time
discrete dynamical systems over a finite set $X$ are \textit{%
non-deterministic finite state automata} (e.g., \cite{FiniteStateAut}) and 
\textit{sequential dynamical systems} \cite{SDS}.

\textit{Deterministic} time discrete dynamical systems over a \textit{finite
field} are mappings $f:K^{n}\rightarrow K^{n},$ where $K$ is a finite field
and $n\in 
\mathbb{N}
$ the dimension of the system. They constitute a particular class of
deterministic time discrete dynamical systems over a finite set $X$, namely,
the class in which the finite set $X$ can be endowed with the algebraic
structure of a finite field. This property allows for a richer mathematical
framework within which these systems can be studied. For instance, it can be
shown that every component function $f_{i}:K^{n}\rightarrow K$ is a
polynomial function of bounded degree in $n$ variables (see, for example,
pages 368-369 in \cite{MR1429394} or 3.1 in \cite{MR???????}).

The study of dynamical systems generally addresses the question of the
system's long term behavior, in particular, the existence of \textit{fixed
points} and \textit{(limit) cyclic trajectories}. (The state of the system
evolves by iteration of the function $f$ starting from given initial
conditions $x_{0}\in K^{n}.$) In this paper we provide an algebraic and
graph theoretic framework to study a very specific class of nonlinear time
discrete dynamical systems over a finite field, namely, \textit{monomial
dynamical systems over a finite field}. In such systems, every component
function $f_{i}:K^{n}\rightarrow K$ is a \textit{monic nonzero monomial
function}.

Some types of monomial systems and their dynamic behavior have been studied
before: monomial cellular automata \cite{CA}, \cite{MCA}, Boolean monomial
systems \cite{MR2112694}, monomial systems over the p-adic numbers \cite%
{p-adic}, \cite{Fuzzy} and monomial systems over a finite field \cite{Vasiga}%
, \cite{ThesisOmar}, \cite{MR2293353}. \cite{MR2112694} proved a necessary
and sufficient condition for Boolean monomial systems to be fixed point
systems (systems in which all trajectories end in steady states)\footnote{%
This problem is refered to as the \textit{steady state system problem}, see 
\cite{JUST}.}. This condition could be algorithmically exploited. Indeed,
the authors make some suggestive comments in that direction (see 4.3 in \cite%
{MR2112694}). Moreover, the paper describes the structure of the limit
cycles of a special type of Boolean monomial systems. \cite{MR2293353}
presents a necessary and sufficient condition for monomial systems over a
finite field to be fixed point systems. However, this condition is not
easily verifiable and therefore the theorem does not yield a tractable
algorithm in a straightforward way.

Our work was strongly influenced by \cite{MR2112694}, \cite{MR2293353} and 
\cite{ThesisOmar}. However, we took a slightly different approach. The
mathematical formalism we developed allows for a deeper understanding of
monomial dynamical systems over a finite field. In particular, we present an
algorithm of \textit{polynomial complexity} to test whether a given monomial
dynamical system over a finite field is a fixed point system. Furthermore,
we obtain additional theorems that complement the work of \cite{MR2112694}, 
\cite{MR2293353} and provide alternative proofs to many results in \cite%
{MR2112694}. Our formalism also constitutes a basis for the study of \textit{%
monomial control systems}, to be presented elsewhere.

It is pertinent to mention the work of \cite{Elspas} regarding \textit{linear%
} time discrete dynamical systems over a finite field, in which the number
of limit cycles and their lengths is linked to the factorization (in so
called elementary divisor polynomials) of the characteristic polynomial of
the matrix representing the system. (See also \cite{MR2175374} for a more
mathematical exposition and \cite{FiniteStateAut}, \cite{REGER&SCHMIDT} for
applications of the Boolean case in control theory.) \ Furthermore, in \cite%
{AffineSys}, the affine case (a linear map followed by a translation) was
studied. An interesting contribution was made by Paul Cull (\cite{CULL}),
who extended the considerations to nonlinear functions, and showed how to
reduce them to the linear case. However, Cull's approach does not yield an
algorithm of polynomial complexity to solve the steady state system problem.
Moreover, according to \cite{JUST}, this might in general not be possible as
a matter of principle.

The organization of this article is the following:

Section 2 establishes an algebraic and graph theoretic framework within
which monomial dynamical systems over a finite field are studied. It starts
with some basic definitions and algebraic results (some of which are proved
in the appendix) and leads the reader to the first important result: Theorem
2, which states that the monoid of $n$-dimensional monomial dynamical
systems over a finite field $K$ is isomorphic to a certain monoid of
matrices. Section 2 finishes with propositions about the relationship
between the matrix $F$ corresponding to a monomial system $f$ (via the
isomorphism mentioned above) and the \textit{adjacency matrix of the
dependency graph} of $f$ (to be defined below).

Section 3 is devoted to the characterization of fixed point systems. These
characterizations are stated in terms of \textit{connectedness properties}
of the dependency graph. We provide some necessary and sufficient conditions
for a system to be a fixed point system (Theorems 6 and 8). Moreover, we
prove several sufficient conditions for special classes of monomial
dynamical systems over a finite field $K.$

Section 4 presents an algorithm of polynomial complexity to test whether a
given monomial dynamical system over a finite field $K$ is a fixed point
system. A detailed complexity analysis of the algorithm is provided.

\section{Algebraic and graph theoretic formalism}

In this section we will introduce the monoid of $n$-dimensional monomial
dynamical systems over a finite field $\mathbf{F}_{q}.$ Furthermore we will
show that this monoid is isomorphic to a certain monoid of matrices. This
result establishes that the composition $f\circ g$ of two monomial dynamical
systems $f,g$ is completely captured by the product $F\cdot G$ of their
corresponding matrices. In addition, we will introduce the concept of
dependency graph of a monomial dynamical system $f$ and prove that the
adjacency matrix of the dependency graph is precisely the matrix $F$
associated with $f$ via the isomorphism mentioned above. This finding allows
us to link topological properties of the dependency graph with the dynamics
of $f$.

\begin{definition}[Notational Definition]
Since for every finite field $K$ there is a prime number $p\in 
\mathbb{N}
$ (the characteristic of $K$) and a natural number $n\in 
\mathbb{N}
$ such that for the number of elements $\left\vert K\right\vert $ of $K$ it
holds%
\begin{equation*}
\left\vert K\right\vert =p^{n}
\end{equation*}%
we will denote a finite field with $\mathbf{F}_{q}$, where $q$ stands for
the number of elements of the field. It is of course understood that $q$ is
a power of the (prime) characteristic of the field.
\end{definition}

\begin{definition}
\label{Def.ExpSet}Let $\mathbf{F}_{q}$ be a finite field. The set%
\begin{equation*}
E_{q}:=\{0,...q-1\}\subset 
\mathbb{N}%
_{0}
\end{equation*}%
is called the exponents set to the field $\mathbf{F}_{q}.$
\end{definition}

\begin{definition}
\label{Def.MonDynSyst}Let $\mathbf{F}_{q}$ be a finite field. A map $f:%
\mathbf{F}_{q}^{n}\rightarrow \mathbf{F}_{q}^{n}$ is called a monomial
dynamical system over $\mathbf{F}_{q}$ if for every $i\in \{1,...,n\}$ there
exists a tuple $(F_{i1},...,F_{in})\in E_{q}^{n}$ such that%
\begin{equation*}
f_{i}(x)=x_{1}^{F_{i1}}...x_{n}^{F_{in}}\text{ }\forall \text{ }x\in \mathbf{%
F}_{q}^{n}
\end{equation*}
\end{definition}

\begin{remark}
\label{Rem.ZeroExcludedFromMonSyst.}As opposed to \cite{MR2112694}, we
exclude in the definition of monomial dynamical system the possibility that
one of the functions $f_{i}$ is equal to the zero function. However, in
contrast to \cite{MR2293353}, we do allow the case $f_{i}\equiv 1$ in our
definition. This is not a loss of generality because of the following: If we
were studying a dynamical system $f:\mathbf{F}_{q}^{n}\rightarrow \mathbf{F}%
_{q}^{n}$ where one of the functions, say $f_{j}$, was equal to zero then
for every initial state $x\in \mathbf{F}_{q}^{n}$ after one iteration the
system would be in a state $f(x)$ whose $j$th entry is zero. In all
subsequent iterations the value of the $j$th entry would remain zero. As a
consequence, the long term dynamics of the system are reflected in the
projection%
\begin{equation*}
\pi _{\hat{\jmath}}(y):=(y_{1},...,y_{j-1},y_{j+1},...,y_{n})^{t}
\end{equation*}%
and it is sufficient to study the system%
\begin{eqnarray*}
\widetilde{f} &:&\mathbf{F}_{q}^{n-1}\rightarrow \mathbf{F}_{q}^{n-1} \\
y &\mapsto &%
\begin{pmatrix}
f_{1}(y_{1},...,y_{j-1},0,y_{j+1},...,y_{n}) \\ 
\vdots  \\ 
f_{j-1}(y_{1},...,y_{j-1},0,y_{j+1},...,y_{n}) \\ 
f_{j+1}(y_{1},...,y_{j-1},0,y_{j+1},...,y_{n}) \\ 
\vdots  \\ 
f_{n}(y_{1},...,y_{j-1},0,y_{j+1},...,y_{n})%
\end{pmatrix}%
\end{eqnarray*}%
In general, this system $\widetilde{f}$ could contain component functions
equal to the zero function, since every component $f_{i}$ that depends on
the variable $x_{j}$ would become zero. As a consequence, the procedure
described above needs to be applied several times until the lower $n^{\prime
}$-dimensional system obtained does not contain component functions equal to
zero. The long term dynamics of $f$ are reflected in the projection to an $%
n^{\prime }$-dimensional subspace, in particular, all the cycles and fixed
points of $f$ are located in this lower dimensional space. Moreover, points
located outside this lower dimensional subspace are transient states of the
system. It is also possible that this repeated procedure yields the one
dimensional zero function. In this case, we can conclude that the original
system $f$ is a fixed point system with $(0,...,0)\in \mathbf{F}_{q}^{n}$ as
its unique fixed point. Note that this procedure reduces the dimension by $s$
where $0\leq s\leq n.$ As a consequence, the procedure needs to be iterated
at most $n$ times.
\end{remark}

As stated in Theorem 16 and Theorem 20 of \cite{MR???????}, every function $%
h:\mathbf{F}_{q}^{n}\rightarrow \mathbf{F}_{q}$ is a polynomial function in $%
n$ variables where no variable appears to a power higher or equal to $q.$
Calculating the composition of a dynamical system $f:\mathbf{F}%
_{q}^{n}\rightarrow \mathbf{F}_{q}^{n}$ with itself, we face the situation
where some of the exponents exceed the value $q-1$ and need to be reduced
according to the well-known rule%
\begin{equation}
a^{q}=a\text{ }\forall \text{ }a\in \mathbf{F}_{q}  \label{eq.ExponentRule}
\end{equation}%
This process can be accomplished systematically if we look at the power $%
x^{p}$ (where $p>q$) as a polynomial in the ring $\mathbf{F}_{q}[\tau ]$ as
described in the Lemma and Definition below. But first we need an auxiliary
result:

\begin{lemma}
\label{ExponentLemma}Let $\mathbf{F}_{q}$ be a finite field and $a\in 
\mathbb{N}
_{0}$ a nonnegative integer. Then%
\begin{equation*}
x^{a}=1\text{ }\forall \text{ }x\in \mathbf{F}_{q}\backslash
\{0\}\Leftrightarrow \exists \text{ }\lambda \in 
\mathbb{N}
_{0}:a=\lambda (q-1)
\end{equation*}
\end{lemma}

\begin{proof}
If $a=\lambda (q-1)$ then $x^{a}=x^{\lambda (q-1)}=(x^{(q-1)})^{\lambda }=1$ 
$\forall $ $x\in \mathbf{F}_{q}\backslash \{0\}$ by (\ref{eq.ExponentRule}).
Now assume $x^{a}=1$ $\forall $ $x\in \mathbf{F}_{q}\backslash \{0\}$ and
write $a=\alpha (q-1)+s$ with suitable $\alpha \in 
\mathbb{N}
_{0}$ and $0\leq s\leq (q-1).$ Then it follows%
\begin{equation*}
1=x^{a}=x^{\lambda (q-1)+s}=x^{\lambda (q-1)}x^{s}=x^{s}\text{ }\forall 
\text{ }x\in \mathbf{F}_{q}\backslash \{0\}
\end{equation*}%
As a consequence, the polynomial $\tau ^{s}-\tau ^{0}\in \mathbf{F}_{q}[\tau
]$ has%
\begin{equation*}
\left\vert \mathbf{F}_{q}\right\vert -1=q-1\geq s=\deg (\tau ^{s}-\tau )
\end{equation*}%
roots in $\mathbf{F}_{q}$ and must be therefore of degree $s=q-1.$
Thus  $a=(\alpha +1)(q-1).$
\end{proof}

\begin{lemma}[and Definition]
\label{red.alg.Lemma}Let $\mathbf{F}_{q}$ be a finite field and $c\in 
\mathbb{N}
_{0}$ a nonnegative integer. The degree of the (unique) remainder of the
polynomial division $\tau ^{c}\div (\tau ^{q}-\tau )$ is called $red_{q}(c).$
$red_{q}(c)$ satisfies the following properties

\begin{enumerate}
\item $red_{q}(red_{q}(c))=red_{q}(c)$

\item $red_{q}(c)=0$ $\Leftrightarrow c=0$

\item For $a,b\in 
\mathbb{N}
_{0},$ $x^{a}=x^{b}$ $\forall $ $x\in \mathbf{F}_{q}\Leftrightarrow
red_{q}(a)=red_{q}(b)$

\item For $a,b\in 
\mathbb{N}
,$ $red_{q}(a)=red_{q}(b)\Leftrightarrow \exists $ $\alpha \in 
\mathbb{Z}
:a=b+\alpha (q-1)$
\end{enumerate}
\end{lemma}

\begin{proof}
By the division algorithm there are unique $g,r\in \mathbf{F}_{q}[\tau ]$
with either $r=0$ or $\deg (r)<\deg (\tau ^{q}-\tau )$ such that%
\begin{equation*}
\tau ^{c}=g(\tau ^{q}-\tau )+r
\end{equation*}%
If we look at the corresponding polynomial functions\footnote{%
If $r\in \mathbf{F}_{q}[\tau ]$ is a polynomial of degree $n,$ i.e. $%
r=\sum\limits_{i=0}^{n}a_{i}\tau ^{i},$ then $\widetilde{r}$ is defined as
the polynomial function%
\begin{eqnarray*}
\widetilde{r} &:&\mathbf{F}_{q}\rightarrow \mathbf{F}_{q} \\
x &\mapsto &\sum\limits_{i=0}^{n}a_{i}x^{i}
\end{eqnarray*}%
} defined on $\mathbf{F}_{q}$ it follows by (\ref{eq.ExponentRule})%
\begin{equation}
x^{c}=\widetilde{r}(x)\text{ }\forall \text{ }x\in \mathbf{F}_{q}
\label{eq. rest}
\end{equation}%
In particular, $r\neq 0.$ From the division process it is also clear that $r$
must be a monomial and we conclude $r=\tau ^{red_{q}(c)}$ with $%
red_{q}(c)<q. $ The first property follows trivially from the fact $%
red_{q}(c)<q.$ The second property follows immediately from evaluating the
equation $x^{c}=x^{red_{q}(c)}$ (i.e. equation (\ref{eq. rest})) at the
value $x=0.$ The third property is shown as follows: By the division
algorithm $\exists _{1}$ $g_{a}$ $,\ g_{b}$ $,$ $r_{a}$ $,\ r_{b}$ $\in 
\mathbf{F}_{q}[\tau ]$ such that%
\begin{eqnarray}
\tau ^{a} &=&g_{a}(\tau ^{q}-\tau )+r_{a}=g_{a}(\tau ^{q}-\tau )+\tau
^{red_{q}(a)}  \label{poly.eq.} \\
\tau ^{b} &=&g_{b}(\tau ^{q}-\tau )+r_{b}=g_{b}(\tau ^{q}-\tau )+\tau
^{red_{q}(b)}  \notag
\end{eqnarray}%
From $x^{a}=x^{b}$ $\forall $ $x\in \mathbf{F}_{q}$ now we have%
\begin{equation*}
x^{red_{q}(a)}=x^{red_{q}(b)}\text{ }\forall x\in \mathbf{F}_{q}
\end{equation*}%
and since $red_{q}(a),red_{q}(b)<q$ we get $red_{q}(a)=red_{q}(b).$ On the
other hand, from $red_{q}(a)=red_{q}(b)$ it would follow from equations (\ref%
{poly.eq.})%
\begin{equation*}
\tau ^{a}-g_{a}(\tau ^{q}-\tau )=\tau ^{b}-g_{b}(\tau ^{q}-\tau )
\end{equation*}%
and thus by (\ref{eq.ExponentRule})%
\begin{equation*}
x^{a}=x^{b}\text{ }\forall \text{ }x\in \mathbf{F}_{q}
\end{equation*}%
Last we prove the fourth claim: If $red_{q}(a)=red_{q}(b)$ then by 3. we have%
\begin{equation*}
x^{a}=x^{b}\text{ }\forall \text{ }x\in \mathbf{F}_{q}
\end{equation*}%
Now assume wlog $a\geq b$ and $d:=a-b\in 
\mathbb{N}
_{0}.$ Then the last equation can be written as%
\begin{equation*}
x^{b}x^{d}=x^{b}\text{ }\forall \text{ }x\in \mathbf{F}_{q}
\end{equation*}%
yielding%
\begin{equation*}
x^{d}=1\text{ }\forall \text{ }x\in \mathbf{F}_{q}\backslash \{0\}
\end{equation*}%
By Lemma \ref{ExponentLemma} we have $\exists $ $\alpha \in 
\mathbb{N}
_{0}:d=\alpha (q-1)$ and therefore  $a=b+\alpha (q-1)$ or $%
b=a-\alpha (q-1).$ Now assume the converse, namely $\exists $ $\alpha \in 
\mathbb{Z}
:a=b+\alpha (q-1).$ Assume wlog $\alpha \geq 0$ (otherwise consider $%
b=a-\alpha (q-1)$). Then we would have%
\begin{equation*}
\tau ^{a}=\tau ^{\alpha (q-1)}\tau ^{b}
\end{equation*}%
and thus by Lemma \ref{ExponentLemma} 
\begin{equation*}
x^{a}=x^{b}\text{ }\forall \text{ }x\in \mathbf{F}_{q}\backslash \{0\}
\end{equation*}%
Since $a,b>0$ we also have%
\begin{equation*}
x^{a}=x^{b}\text{ }\forall \text{ }x\in \mathbf{F}_{q}
\end{equation*}
\end{proof}

\begin{remark}
\label{PowerOfredRemark}From the properties above we have $x^{a}=$ $%
x^{red_{q}(a)}$ $\forall $ $x\in \mathbf{F}_{q}.$
\end{remark}

The "exponents arithmetic" needed when calculating the composition of
dynamical systems $f,g:\mathbf{F}_{q}^{n}\rightarrow \mathbf{F}_{q}^{n}$ can
be formalized based on the reduction algorithm described by the previous
lemma. Indeed, the set%
\begin{equation*}
E_{q}=\{0,1,...,(q-2),(q-1)\}\subset 
\mathbb{Z}%
\end{equation*}%
together with the operations of addition $a\oplus b:=red_{q}(a+b)$ and
multiplication $a\bullet b:=red_{q}(ab)$ is a commutative semiring with
identity $1$. We call this commutative semiring the \textit{exponents
semiring} of the field $\mathbf{F}_{q}.$ This result is proved in the
Appendix (see Theorem \ref{Th.ExponentsSemiring}). We also defer to the
appendix the proof of the following lemma:

\begin{lemma}
Let $n\in 
\mathbb{N}
$ be a natural number, $\mathbf{F}_{q}$ be a finite field and $E_{q}$ the
exponents semiring of $\mathbf{F}_{q}.$ The set $M(n\times n;$ $E_{q})$ of $%
n\times n$ quadratic matrices with entries in the semiring $E_{q}$ together
with the operation $\cdot $ of matrix multiplication (which is defined in
terms of the operations $\oplus $ and $\bullet $ on the matrix entries) over 
$E_{q}$ is a monoid.
\end{lemma}

\begin{remark}[and Definition]
The operation $red_{q}:%
\mathbb{N}
_{0}\rightarrow $ $E_{q}$ can be extended to matrices $M(n\times n;$ $%
\mathbb{N}
_{0})$ by applying $red_{q}$ to the entries of the matrix. We call this
extension $mred_{q}:M(n\times n;$ $%
\mathbb{N}
_{0})\rightarrow M(n\times n;$ $E_{q}).$ See Remark \ref{RemarkMatRed} in
the appendix for further details. One important property of $mred_{q}$ shown
in Remark \ref{RemarkMatRed} is%
\begin{equation}
mred_{q}(A)=0\Leftrightarrow A=0  \label{eq.mred(A)=0Body}
\end{equation}
\end{remark}

\begin{definition}
Let $\mathbf{F}_{q}$ be a finite field and $n,m\in 
\mathbb{N}
$ natural numbers. The set%
\begin{eqnarray*}
MF_{m}^{n}(\mathbf{F}_{q}) &:&=\{f:\mathbf{F}_{q}^{m}\rightarrow \mathbf{F}%
_{q}^{n}\text{ }| \\
\exists \text{ }F &\in &M(n\times
m;E_{q}):f_{i}(x):=x_{1}^{F_{i1}}...x_{m}^{F_{im}}\text{ }\forall \text{ }%
x\in \mathbf{F}_{q}^{n}\}
\end{eqnarray*}%
is called the set of $n$\textit{-dimensional monomial mappings in }$m$%
\textit{\ variables.}
\end{definition}

\begin{lemma}
\label{MonoidIsomorph.Lemma}Let $\mathbf{F}_{q}$ be a finite field and $%
n,m,r\in 
\mathbb{N}
$ natural numbers. Furthermore, let $f\in MF_{n}^{m}(\mathbf{F}_{q})$ and $%
g\in MF_{m}^{r}(\mathbf{F}_{q})$ with%
\begin{eqnarray*}
f_{i}(x) &=&x_{1}^{F_{i1}}...x_{n}^{F_{in}}\text{ }\forall \text{ }x\in 
\mathbf{F}_{q}^{n},\text{ }i=1,...,m \\
g_{j}(x) &=&x_{1}^{G_{j1}}...x_{m}^{G_{jm}}\text{ }\forall \text{ }x\in 
\mathbf{F}_{q}^{n},\text{ }j=1,...,r
\end{eqnarray*}%
where $F\in M(m\times n;E_{q})$ and $G\in M(r\times m;$ $E_{q}).$ Then for
their composition $g\circ f:\mathbf{F}_{q}^{n}\rightarrow \mathbf{F}_{q}^{r}$
it holds%
\begin{equation*}
(g\circ
f)_{k}(x)=\prod\limits_{j=1}^{n}x_{j}{}^{red_{q}(\tsum%
\limits_{l=1}^{m}G_{kl}F_{lj})}\text{ }\forall \text{ }x\in \mathbf{F}%
_{q}^{n},\text{ }k\in \{1,...,r\}
\end{equation*}
\end{lemma}

\begin{proof}
From the definition It follows for every $k\in \{1,...,r\}$%
\begin{equation*}
(g\circ
f)_{k}(x)=\prod\limits_{l=1}^{m}(f_{l}(x))^{G_{kl}}=\prod\limits_{l=1}^{m}(%
\prod\limits_{j=1}^{n}x_{j}^{F_{lj}})^{G_{kl}}
\end{equation*}%
For a fixed but arbitrary $m\in 
\mathbb{N}
$ we will prove the claim using induction on the dimension $n$ of the
mapping $g\circ f$. For $n=1$ we have $(g\circ
f)_{k}(x)=\prod\limits_{l=1}^{m}(x_{1}^{F_{l1}})^{G_{kl}}=\prod%
\limits_{l=1}^{m}x_{1}^{G_{kl}F_{l1}}=x_{1}^{\tsum%
\limits_{l=1}^{m}G_{kl}F_{l1}}=x_{1}^{red_{q}(\tsum%
\limits_{l=1}^{m}G_{kl}F_{l1})}$ (see Remark \ref{PowerOfredRemark}), thus
the claim holds in dimension $1.$ Now we consider the case $n+1:$%
\begin{eqnarray*}
(g\circ f)_{k}(x)
&=&\prod\limits_{l=1}^{m}(\prod\limits_{j=1}^{n+1}x_{j}^{F_{lj}})^{G_{kl}} \\
&=&\prod\limits_{l=1}^{m}(x_{(n+1)}^{F_{l(n+1)}}\prod%
\limits_{j=1}^{n}x_{j}^{F_{lj}})^{G_{kl}} \\
&=&\prod\limits_{l=1}^{m}\left(
x_{(n+1)}^{G_{kl}F_{l(n+1)}}(\prod\limits_{j=1}^{n}x_{j}^{F_{lj}})^{G_{kl}}%
\right) \\
&=&\prod\limits_{l=1}^{m}(x_{(n+1)}^{G_{kl}F_{l(n+1)}})\prod%
\limits_{l=1}^{m}(\prod\limits_{j=1}^{n}x_{j}^{F_{lj}})^{G_{kl}}
\end{eqnarray*}%
and by induction hypothesis%
\begin{eqnarray*}
&=&x_{(n+1)}^{\sum\limits_{l=1}^{m}G_{kl}F_{l(n+1)}}\prod%
\limits_{j=1}^{n}x_{j}{}^{red_{q}(\tsum\limits_{l=1}^{m}G_{kl}F_{lj})} \\
&=&x_{(n+1)}^{\sum\limits_{l=1}^{m}G_{kl}F_{l(n+1)}}\prod%
\limits_{j=1}^{n}x_{j}{}^{\tsum\limits_{l=1}^{m}G_{kl}F_{lj}} \\
&=&\prod\limits_{j=1}^{n+1}x_{j}{}^{\tsum\limits_{l=1}^{m}G_{kl}F_{lj}} \\
&=&\prod\limits_{j=1}^{n+1}x_{j}{}^{red_{q}(\tsum%
\limits_{l=1}^{m}G_{kl}F_{lj})}
\end{eqnarray*}
\end{proof}

\begin{remark}[and Lemma]
\label{Rem.Comp.And.Mat.Mult.}If we generalize the matrix multiplication
defined on the monoid $M(n\times n;$ $E_{q})$ for matrices $F\in M(m\times
n;E_{q})$ and $G\in M(n\times m;$ $E_{q})$ then we can write%
\begin{equation*}
(g\circ f)_{k}(x)=\prod\limits_{j=1}^{n}x_{j}{}^{(G\cdot F)_{kj}}\text{ }%
\forall \text{ }x\in \mathbf{F}_{q}^{n},\text{ }k\in \{1,...,n\}
\end{equation*}%
To see this, apply the Lemmas \ref{AddLemma} and \ref{MonoidIsomorph.Lemma}
as well as the definitions of $\oplus $ and $\bullet $ to $%
\prod\limits_{j=1}^{n}x_{j}{}^{(G\cdot F)_{kj}}:$ 
\begin{eqnarray*}
\prod\limits_{j=1}^{n}x_{j}{}^{(G\cdot F)_{kj}}
&=&\prod\limits_{j=1}^{n}x_{j}{}^{(G_{k1}\bullet F_{1j}\oplus ...\oplus
G_{km}\bullet F_{mj})} \\
&=&\prod\limits_{j=1}^{n}x_{j}{}^{red_{q}(G_{k1}F_{1j})\oplus ...\oplus
red_{q}(G_{km}F_{mj})} \\
&=&\prod%
\limits_{j=1}^{n}x_{j}{}^{red_{q}(red_{q}(G_{k1}F_{1j})+...+red_{q}(G_{km}F_{mj}))}
\\
&=&\prod\limits_{j=1}^{n}x_{j}{}^{red_{q}(\tsum%
\limits_{l=1}^{m}G_{kl}F_{lj})}\text{ } \\
&=&(g\circ f)_{k}(x)
\end{eqnarray*}
\end{remark}

\begin{theorem}
Let $\mathbf{F}_{q}$ be a finite field. The set%
\begin{eqnarray*}
MF_{n}^{n}(\mathbf{F}_{q}) &:&=\{f:\mathbf{F}_{q}^{n}\rightarrow \mathbf{F}%
_{q}^{n}\text{ }| \\
\exists \text{ }F &\in &M(n\times
n;E_{q}):f_{i}(x):=x_{1}^{F_{i1}}...x_{n}^{F_{in}}\text{ }\forall \text{ }%
x\in \mathbf{F}_{q}^{n}\}
\end{eqnarray*}%
of all monomial dynamical systems over $\mathbf{F}_{q}$ together with the
composition $\circ $ of mappings is a monoid.
\end{theorem}

\begin{proof}
By Lemma \ref{MonoidIsomorph.Lemma} the set $MF_{n}^{n}(\mathbf{F}_{q})$ is
closed under composition. Composition of mappings is trivially associative.
The identity function%
\begin{eqnarray*}
Id &:&\mathbf{F}_{q}^{n}\rightarrow \mathbf{F}_{q}^{n} \\
x &\mapsto &x
\end{eqnarray*}%
is a monomial function and is therefore the identity element of the monoid $%
(MF_{n}^{n}(\mathbf{F}_{q}),\circ ).$
\end{proof}

\begin{theorem}
\label{ThmMonoidIsomorphism}The monoids $M(n\times n;$ $E_{q})$ and $%
MF_{n}^{n}(\mathbf{F}_{q})$ are isomorphic.
\end{theorem}

\begin{proof}
From the definition of $MF_{n}^{n}(\mathbf{F}_{q})$ it is clear that the
mapping%
\begin{eqnarray*}
\Psi &:&M(n\times n;E_{q})\rightarrow MF_{n}^{n}(\mathbf{F}_{q}) \\
G &\mapsto &\Psi (G)
\end{eqnarray*}%
such that%
\begin{equation*}
\Psi (G)_{i}(x):=x_{1}^{G_{i1}}...x_{n}^{G_{in}}\text{ for }i=1,...,n
\end{equation*}%
is a bijection. Moreover, $\Psi (I)=id$. In addition, by Remark \ref%
{Rem.Comp.And.Mat.Mult.} it follows easily that%
\begin{equation*}
\Psi (F\cdot G)=\Psi (F)\circ \Psi (G)
\end{equation*}
\end{proof}

\begin{remark}[and Definition]
\label{RemarkMatPower}For a given monomial dynamical system  $f\in
MF_{n}^{n}(\mathbf{F}_{q})$ the matrix $F:=\Psi ^{-1}(f)$ is called the 
\textit{corresponding matrix} of the system $f.$ For a matrix power in the
monoid $M(n\times n;E_{q})$ we use the notation $F^{\cdot m}.$ By induction
it can be easily shown%
\begin{equation*}
\Psi ^{-1}(f^{m})=F^{\cdot m}
\end{equation*}
\end{remark}

\begin{remark}[and Definition]
The image of the $n\times n$ zero matrix  $0\in M(n\times n;E_{q})$
under the isomorphism $\Psi $ has the property%
\begin{equation*}
\Psi (0)(x)_{i}=1\text{ }\forall \text{ }x\in \mathbf{F}_{q}^{n}
\end{equation*}%
we call this monomial function the \textit{one function} $\boldsymbol{1}%
:=\Psi (0).$
\end{remark}

\begin{definition}[Notational Definition]
A directed graph%
\begin{equation*}
G=(V_{G},E_{G},\pi _{G}:E_{G}\rightarrow V_{G}\times V_{G})
\end{equation*}%
that allows self loops and parallel directed edges is called \textit{digraph}%
.
\end{definition}

\begin{definition}
Let $M$ be a nonempty finite set. Furthermore, let $n:=\left\vert
M\right\vert $ be the cardinality of $M.$ A \textit{numeration} of the
elements of $M$ is a bijective mapping%
\begin{equation*}
f:M\rightarrow \{1,...,n\}
\end{equation*}%
Given a numeration $f$ of the set $M$ we write%
\begin{equation*}
M=\{f_{1},...,f_{n}\}
\end{equation*}%
where the unique element $x\in M$ with the property $f(x)=i\in $ $%
\{1,...,n\} $ is denoted as $f_{i}$.
\end{definition}

\begin{definition}[Notational Definition]
\label{Def.mon.Dep.Graph}Let $f\in MF_{n}^{n}(\mathbf{F}_{q})$ be a monomial
dynamical system and $G=(V_{G},$ $E_{G},$ $\pi _{G})$ a digraph with vertex
set $V_{G}$ of cardinality $\left\vert V_{G}\right\vert =n.$ Furthermore,
let $F:=\Psi ^{-1}(f)$ be the corresponding matrix of $f.$ The digraph $G$
is called \emph{dependency graph} of $f$ iff a numeration $a:M\rightarrow
\{1,...,n\}$ of the elements of $V_{G}$ exists such that $\forall $ $i,j\in
\{1,...,n\}$ there are \textbf{exactly} $F_{ij}$ directed edges $%
a_{i}\rightarrow a_{j}$ in the set $E_{G},$ i.e.%
\begin{equation*}
\left\vert \pi _{f}^{-1}((a_{i},a_{j}))\right\vert =F_{ij}
\end{equation*}
\end{definition}

\begin{remark}
It is easy to show that if $G$ and $H$ are dependency graphs of $f$ then $G$
and $H$ are isomorphic. In this sense we speak from \textit{the} dependency
graph of $f$ and denote it by $G_{f}=(V_{f},$ $E_{f},$ $\pi _{f}).$ Our
definition of dependency graph differs slightly from the definition used in 
\cite{MR2112694}.
\end{remark}

\begin{definition}[Notational Definition]
\label{Def.Sequence}Let $G=(V_{G},$ $E_{G},$ $\pi _{G})$ be a digraph. Two
vertices $a,b\in V_{G}$ are called \textit{connected} if there is a $t\in 
\mathbb{N}
_{0}$ and (not necessarily different) vertices $v_{1},...,v_{t}\in V_{G}$
such that%
\begin{equation*}
a\rightarrow v_{1}\rightarrow v_{2}\rightarrow ...\rightarrow
v_{t}\rightarrow b
\end{equation*}%
In this situation we write $a\rightsquigarrow _{s}b,$ where $s$ is the
number of directed edges involved in the \textit{sequence} from $a$ to $b$
(in this case $s=t+1$). Two sequences $a\rightsquigarrow _{s}b$ of the same
length are considered different if the directed \textit{edges} involved are
different or the order at which they appear is different, even if the
visited vertices are the same. As a convention, a single vertex $a\in V_{G}$
is always connected to itself $a\rightsquigarrow _{0}a$ by an empty sequence
of length $0.$
\end{definition}

\begin{definition}[Notational Definition]
\label{Def.Path}Let $G=(V_{G},$ $E_{G},$ $\pi _{G})$ be a digraph and $%
a,b\in V_{G}$ two vertices. A sequence $a\rightsquigarrow _{s}b$%
\begin{equation*}
a\rightarrow v_{1}\rightarrow v_{2}\rightarrow ...\rightarrow
v_{t}\rightarrow b
\end{equation*}%
is called a \emph{path}, if no vertex $v_{i}$ is visited more than once. If $%
a=b,$ but no other vertex is visited more than once, $a\rightsquigarrow
_{s}b $ is called a \emph{closed path}.
\end{definition}

\begin{definition}
Let $G=(V_{G},$ $E_{G},$ $\pi _{G})$ be a digraph. Two vertices  $%
a,b\in V_{G}$ are called \textit{strongly connected }if there are natural
numbers $s,t\in 
\mathbb{N}
$ such that%
\begin{equation*}
a\rightsquigarrow _{s}b\text{ and }b\rightsquigarrow _{t}a
\end{equation*}%
In this situation we write $a\rightleftharpoons b.$
\end{definition}

\begin{theorem}[and Definition]
Let $G=(V_{G},$ $E_{G},$ $\pi _{G})$ be a digraph. $\rightleftharpoons $ is
an equivalence relation on $V_{G}$ called \textit{strong equivalence}. The
equivalence class of any vertex $a\in V_{G}$ is called a \textit{strongly
connected component} and denoted by $\overleftrightarrow{a}\subseteq V_{G}.$
\end{theorem}

\begin{proof}
See Definition 3.1 (2) in \cite{MR2112694}.
\end{proof}

\begin{definition}
Let $G=(V_{G},$ $E_{G},$ $\pi _{G})$ be a digraph and $a\in V_{G}$ one of
its vertices. The \textit{strongly connected component }$\overleftrightarrow{%
a}\subseteq V_{G}$ is called trivial iff $\overleftrightarrow{a}=\{a\}$ and
there is no edge $a\rightarrow a$ in $E_{G}$.
\end{definition}

\begin{definition}
Let $G=(V_{G},$ $E_{G},$ $\pi _{G})$ be a digraph with vertex set $V_{G}$ of
cardinality $\left\vert V_{G}\right\vert =n$ and $V_{G}=\{a_{1},...,a_{n}\}$
a numeration of the elements of $V_{G}$. The matrix $A\in M(n\times n;$ $%
\mathbb{N}
_{0})$ whose entries are defined as%
\begin{equation*}
A_{ij}:=\text{number of edges }a_{i}\rightarrow a_{j}\text{ contained in }%
E_{G}
\end{equation*}%
for $i,j=1,...,n$ is called \textit{adjacency matrix} of $G$ with the
numeration $a.$
\end{definition}

\begin{theorem}
\label{ThmPowerOfAdjMat}Let $G=(V_{G},$ $E_{G},$ $\pi _{G})$ be a digraph
with vertex set $V_{G}$ of cardinality $\left\vert V_{G}\right\vert =n$ and $%
V_{G}=\{a_{1},...,a_{n}\}$ a numeration of the elements of $V_{G}$.
Furthermore, let $A\in M(n\times n;$ $%
\mathbb{N}
_{0})$ be its adjacency matrix (with the numeration $a$), $m\in 
\mathbb{N}
$ a natural number and%
\begin{equation*}
B:=A^{m}\in M(n\times n;%
\mathbb{N}
_{0})
\end{equation*}%
the $m$th power of $A.$ Then $\forall $ $i,j\in \{1,...,n\}$ the entry $%
B_{ij}$ of $B$ is equal to the number of different sequences $%
a_{i}\rightsquigarrow _{m}a_{j}$ of length $m.$
\end{theorem}

\begin{proof}
The proof of this well-known result can be found in \cite{MR0184874}.
\end{proof}

\begin{remark}
\label{RemarkPowersOfF}Let $f\in MF_{n}^{n}(\mathbf{F}_{q})$ be a monomial
dynamical system. Furthermore, let $G_{f}=(V_{f},$ $E_{f},$ $\pi _{f})$ the
dependency graph of $f$ and $V_{f}=\{a_{1},...,a_{n}\}$ the associated
numeration of the elements of $V_{f}.$ Then, according to the definition of
dependency graph, $F:=\Psi ^{-1}(f)$ (the corresponding matrix of $f$) is
precisely the adjacency matrix of $G_{f}$ with the numeration $a.$ Now, by
Remarks \ref{RemarkMatPower} and \ref{RemarkMatRed} we can conclude%
\begin{equation}
\Psi ^{-1}(f^{m})=mred_{q}(F^{m})  \label{eq.PowersOfF}
\end{equation}
\end{remark}

\section{Characterization of fixed point systems}

The results proved in the previous section allow us to link topological
properties of the dependency graph with the dynamics of $f$. We will exploit
this feature in this subsection to prove some characterizations of fixed
point systems stated in terms of connectedness properties of the dependency
graph. At the end of this section we also provide a more algebraic
sufficient condition.

\begin{theorem}
\label{Th.FPSwhenOnlyTrivialSCC}Let $\mathbf{F}_{q}$ be a finite field and $%
f\in MF_{n}^{n}(\mathbf{F}_{q})$ a monomial dynamical system. Then $f$ is a
fixed point system with $(1,...,1)^{t}\in \mathbf{F}_{q}^{n}$ as its only
fixed point if and only if its dependency graph only contains trivial
strongly connected components .
\end{theorem}

\begin{proof}
By Remark \ref{RemarkPowersOfF}, $F:=\Psi ^{-1}(f)$ is the adjacency matrix
of the dependency graph of $f.$ If the dependency graph does not contain any
nontrivial strongly connected components, every sequence $a\rightsquigarrow
_{s}b$ between two arbitrary vertices can be at most of length $n-1.$ (A
sequence that revisits a vertex would contain a closed sequence, which is
strongly connected.) Therefore, by theorem \ref{ThmPowerOfAdjMat} $\exists $ 
$m\in 
\mathbb{N}
$ with $m\leq n$ such that $F^{m}=0$ (the zero matrix in $M(n\times n;$ $%
\mathbb{N}
_{0})$). Now, according to equation (\ref{eq.PowersOfF}) we have%
\begin{equation*}
\Psi ^{-1}(f^{m})=mred_{q}(F^{m})=mred_{q}(0)=0
\end{equation*}%
and consequently%
\begin{equation*}
\Psi ^{-1}(f^{r})=0\text{ }\forall \text{ }r\geq m
\end{equation*}%
Thus%
\begin{equation*}
f^{r}=\boldsymbol{1}\text{ }\forall \text{ }r\geq m
\end{equation*}%
If, on the other hand, there is an $m\in 
\mathbb{N}
$ such that%
\begin{equation*}
f^{m+\lambda }=f^{m}=\boldsymbol{1}\text{ }\forall \text{ }\lambda \in 
\mathbb{N}%
\end{equation*}%
applying the isomorphism $\Psi ^{-1}$ (see Remark \ref{RemarkMatPower}) we
obtain%
\begin{equation*}
F^{\cdot (m+\lambda )}=F^{\cdot m}=0\text{ }\forall \text{ }\lambda \in 
\mathbb{N}%
\end{equation*}%
and (see equation (\ref{eq.PowersOfF}))%
\begin{equation*}
mred_{q}(F^{m+\lambda })=mred_{q}(F^{m})=0\text{ }\forall \text{ }\lambda
\in 
\mathbb{N}%
\end{equation*}%
It follows from equation (\ref{eq.mred(A)=0Body}) (See also Remark \ref%
{RemarkMatRed})%
\begin{equation*}
F^{m+\alpha }=0\text{ }\forall \text{ }\alpha \in 
\mathbb{N}
_{0}
\end{equation*}%
Now by theorem \ref{ThmPowerOfAdjMat} there are no sequences $%
a\rightsquigarrow _{s}b$ between any two arbitrary vertices $a,b$ of length
larger than $m-1.$ As a consequence, there cannot be any nontrivial strongly
connected components in the dependency graph of $f.$
\end{proof}

\begin{definition}
A monomial dynamical system $f\in MF_{n}^{n}(\mathbf{F}_{q})$ whose
dependency graph\textit{\ }contains nontrivial strongly connected components
is called \textit{coupled monomial dynamical system}.
\end{definition}

\begin{definition}
Let $G=(V_{G},$ $E_{G},$ $\pi _{G})$ be a digraph, $m\in 
\mathbb{N}
$ a natural number and $a,b\in V_{G}$ two vertices. The number of different
sequences of length $m$ from $a$ to $b$ is denoted by $s_{m}(a,b)\in 
\mathbb{N}
_{0}$ $.$
\end{definition}

\begin{remark}
Let $G=(V_{G},$ $E_{G},$ $\pi _{G})$ be a digraph with vertex set $V_{G}$ of
cardinality $n:=\left\vert V_{G}\right\vert $ and $V_{G}=\{a_{1},...,a_{n}\}$
a numeration of the elements of $V_{G}$. Furthermore, let $m\in 
\mathbb{N}
$ be a natural number and $A\in M(n\times n;$ $%
\mathbb{N}
_{0})$ the adjacency matrix of $G$ with the numeration $a$. Then by Theorem %
\ref{ThmPowerOfAdjMat} we have%
\begin{equation*}
s_{m}(a_{i},a_{j})=(A^{m})_{ij}
\end{equation*}
\end{remark}

\begin{theorem}
Let $\mathbf{F}_{q}$ be a finite field, $f\in MF_{n}^{n}(\mathbf{F}_{q})$ a
coupled monomial dynamical system and $G_{f}=(V_{f},$ $E_{f},$ $\pi _{f})$
its dependency graph. Then $f$ is a fixed point system if and only if there
is an $m\in 
\mathbb{N}
$ such that the following two conditions hold

\begin{enumerate}
\item {\normalsize For every pair of nodes }$a,b${\normalsize \ }$\in V_{f}$%
{\normalsize \ with }$a\rightsquigarrow _{m}b${\normalsize \ there exists
for every }$\lambda \in 
\mathbb{N}
${\normalsize \ an }$a_{\lambda }\in 
\mathbb{Z}
${\normalsize \ such that }$s_{m+\lambda }(a,b)=s_{m}(a,b)+a_{\lambda
}(q-1)\neq 0.$

\item {\normalsize For every pair of nodes }$a,b${\normalsize \ }$\in V_{f}$%
{\normalsize \ with }$s_{m}(a,b)=0${\normalsize \ it holds }$s_{m+\lambda
}(a,b)=0${\normalsize \ }$\forall ${\normalsize \ }$\lambda \in 
\mathbb{N}
.$
\end{enumerate}
\end{theorem}

\begin{proof}
Let $V_{f}=\{a_{1},...,a_{n}\}$ be the numeration of the vertices. If $f$ is
a fixed point system, $\exists $ $m\in 
\mathbb{N}
$ such that%
\begin{equation*}
f^{m+\lambda }=f^{m}\text{ }\forall \text{ }\lambda \in 
\mathbb{N}%
\end{equation*}%
By applying the homomorphism $\Psi ^{-1}$ we get (see Remark \ref%
{RemarkMatPower})%
\begin{equation}
F^{\cdot (m+\lambda )}=F^{\cdot m}\text{ }\forall \text{ }\lambda \in 
\mathbb{N}
\label{eq.Fodes}
\end{equation}%
By Remark \ref{RemarkPowersOfF} it follows%
\begin{equation*}
mred_{q}(F^{m+\lambda })=mred_{q}(F^{m})\text{ }\forall \text{ }\lambda \in 
\mathbb{N}%
\end{equation*}%
Let $i,j\in \{1,...,n\}.$ If, on the one hand, $(F^{\cdot m})_{ij}=0$ then
by (\ref{eq.Fodes}) we would have $(F^{\cdot (m+\lambda )})_{ij}=0$ $\forall 
$ $\lambda \in 
\mathbb{N}
.$ Consequently, by 2. of Lemma \ref{red.alg.Lemma} we have%
\begin{equation*}
(F^{m+\alpha })_{ij}=0\text{ }\forall \text{ }\alpha \in 
\mathbb{N}
_{0}
\end{equation*}%
Now by theorem \ref{ThmPowerOfAdjMat} there are no sequences $%
a_{i}\rightsquigarrow _{s}a_{j}$ of length larger than $m-1.$ In other
words, 2. follows. If, on the other hand, $(F^{\cdot m})_{ij}\neq 0$ then by
(\ref{eq.Fodes}) we would have $(F^{\cdot (m+\lambda )})_{ij}=(F^{\cdot
m})_{ij}\neq 0$ $\forall $ $\lambda \in 
\mathbb{N}
.$ Consequently, by 2. and 4. of Lemma \ref{red.alg.Lemma} $\exists $ $%
a_{\lambda }\in 
\mathbb{Z}
$ such that%
\begin{equation*}
(F^{m+\lambda })_{ij}=(F^{m})_{ij}+a_{\lambda }(q-1)\text{ }\forall \text{ }%
\lambda \in 
\mathbb{N}%
\end{equation*}%
In other words 1. follows. To show the converse we start from the following
fact: Given 1. and 2. and according to Theorem \ref{ThmPowerOfAdjMat} and
Remark \ref{RemarkPowersOfF}%
\begin{equation*}
\text{If }(F^{m})_{ij}=0,\text{ then }(F^{m+\lambda })_{ij}=(F^{m})_{ij}%
\text{ }\forall \text{ }\lambda \in 
\mathbb{N}%
\end{equation*}%
and%
\begin{equation*}
\text{if }(F^{m})_{ij}\neq 0,\text{ then }\exists \text{ }a_{\lambda }\in 
\mathbb{Z}
:(F^{m+\lambda })_{ij}=(F^{m})_{ij}+a_{\lambda }(q-1)\neq 0\text{ }\forall 
\text{ }\lambda \in 
\mathbb{N}%
\end{equation*}%
Now by 2. and 4. of Lemma \ref{red.alg.Lemma} we have%
\begin{equation*}
mred_{q}(F^{m+\lambda })=mred_{q}(F^{m})\text{ }\forall \text{ }\lambda \in 
\mathbb{N}%
\end{equation*}%
and by \ref{RemarkPowersOfF}%
\begin{equation*}
F^{\cdot (m+\lambda )}=F^{\cdot m}\text{ }\forall \text{ }\lambda \in 
\mathbb{N}%
\end{equation*}%
Thus, after applying the isomorphism $\Psi $%
\begin{equation*}
f^{m+\lambda }=f^{m}\text{ }\forall \text{ }\lambda \in 
\mathbb{N}%
\end{equation*}
\end{proof}

The following parameter for digraphs was introduced by \cite{MR2112694}:

\begin{definition}
\label{DefLoopNumber}Let $G=(V_{G},$ $E_{G},$ $\pi _{G})$ be a digraph and $%
a\in V_{G}$ one of its vertices. The number%
\begin{equation*}
\tciLaplace (a):=\min_{\substack{ a\rightsquigarrow _{u}a  \\ %
a\rightsquigarrow _{v}a  \\ u\neq v}}\left\vert u-v\right\vert
\end{equation*}%
is called the loop number of $a.$ If there is no sequence of positive length
from $a$ to $a,$ then $\tciLaplace (a)$ is set to zero.
\end{definition}

\begin{lemma}[and Definition]
Let $G=(V_{G},$ $E_{G},$ $\pi _{G})$ be a digraph and $a\in V_{G}$ one of
its vertices. If $\overleftrightarrow{a}$ is nontrivial then for every $b\in 
\overleftrightarrow{a}$ it holds%
\begin{equation*}
\tciLaplace (b)=\tciLaplace (a)
\end{equation*}%
Therefore, we introduce the loop number of strongly connected components as%
\begin{equation*}
\tciLaplace (\overleftrightarrow{a}):=\tciLaplace (a)
\end{equation*}
\end{lemma}

\begin{proof}
See the proof of Lemma 4.2 in \cite{MR2112694}.
\end{proof}

\begin{remark}
The loop number of any trivial strongly connected component is, due to the
convention made in the definition of loop number, equal to zero.
\end{remark}

\begin{corollary}
\label{CorNecessityOfLoopNo1}Let $\mathbf{F}_{q}$ be a finite field, $f\in
MF_{n}^{n}(\mathbf{F}_{q})$ a coupled monomial dynamical system and $%
G_{f}=(V_{f},$ $E_{f},$ $\pi _{f})$ its dependency graph. If $f$ is a fixed
point system then the loop number of each of its nontrivial strongly
connected components is equal to 1.
\end{corollary}

\begin{proof}
Let $m\in 
\mathbb{N}
$ be as in the statement of the previous theorem. Let $\overleftrightarrow{a}%
\subseteq V_{f}$ be a nontrivial strongly connected component. For every $%
b\in \overleftrightarrow{a}$ we have that $b$ is strongly connected with
itself. Therefore, for every $s\in 
\mathbb{N}
$ there is a $t\geq s$ such that $b\rightsquigarrow _{t}b.$ In particular,
there must be a $u\in 
\mathbb{N}
$ with $u>m$ such that $b\rightsquigarrow _{u}b,$ i.e. $s_{u}(b,b)\geq 1.$
By 2. of the previous theorem we know that $s_{m}(b,b)\neq 0,$ otherwise $%
s_{u}(b,b)=0.$ Now from 1. of the previous theorem we know%
\begin{equation*}
\exists \text{ }a_{\lambda }\in 
\mathbb{Z}
:s_{m+\lambda }(b,b)=s_{m}(b,b)+a_{\lambda }(q-1)\neq 0\text{ }\forall \text{
}\lambda \in 
\mathbb{N}%
\end{equation*}%
and in particular%
\begin{equation*}
s_{m+\lambda }(b,b)\neq 0\text{ }\forall \text{ }\lambda \in 
\mathbb{N}%
\end{equation*}%
Therefore, $\forall $ $\lambda \in 
\mathbb{N}
$ there are sequences $b\rightsquigarrow _{m+\lambda }b$. Thus $\tciLaplace (%
\overleftrightarrow{a})=\tciLaplace (b)=1.$
\end{proof}

\begin{definition}
Let $G=(V_{G},$ $E_{G},$ $\pi _{G})$ be a digraph and $a,b\in V_{G}$ two
vertices. The vertex $a$ is called \textit{recurrently connected }to $b$, if
for every $s\in 
\mathbb{N}
$ there is a $u\geq s$ such that $a\rightsquigarrow _{u}b$.
\end{definition}

\begin{lemma}[and Definition]
Let $G=(V_{G},$ $E_{G},$ $\pi _{G})$ be a digraph with vertex set $V_{G}$ of
cardinality $n:=\left\vert V_{G}\right\vert $. Two vertices $a,b\in V_{G}$
are connected through a sequence $a\rightsquigarrow _{t}b$ of length $t>n-1$
if and only if $a$ is recurrently connected to $b.$
\end{lemma}

\begin{proof}
If there is a sequence $a\rightsquigarrow _{t}b$ of length $t>n-1,$ then it
necessarily revisits one of its vertices, in other words, there is a $c\in
V_{G}$ such that%
\begin{equation*}
a\rightsquigarrow _{t}b=a\rightarrow ...\rightarrow c\rightarrow
...\rightarrow c\rightarrow ...\rightarrow b
\end{equation*}%
Now a sequence $a\rightsquigarrow _{t^{\prime }}b$ can be constructed that
repeats the loop around $c$ as many times as desired. The converse follows
immediately from the definition of recurrent connectedness.
\end{proof}

\begin{remark}
\label{RemarkRecurrentConn}Let $G=(V_{G},$ $E_{G},$ $\pi _{G})$ be a digraph
with vertex set $V_{G}$ of cardinality $n:=\left\vert V_{G}\right\vert $.
Then for any two vertices $a,b\in V_{G}$ it holds: Either $a$ is recurrently
connected to $b$ or there is an $m\in 
\mathbb{N}
$ with $m\leq n$ such that no sequence $a\rightsquigarrow _{t}b$ of length $%
t\geq m$ exists.
\end{remark}

\begin{lemma}
\label{Omar'sProposition}Let $G=(V_{G},$ $E_{G},$ $\pi _{G})$ be a digraph
and $U\subseteq V_{G}$ a nontrivial strongly connected component.
Furthermore, let $t:=\tciLaplace (U)$ be the loop number of $U.$ Then for
each $a,b\in U$ there is an $m\in 
\mathbb{N}
$ such that the graph $G$ contains sequences $a\rightsquigarrow _{m+\lambda
t}b$ of length $m+\lambda t$ $\forall $ $\lambda \in 
\mathbb{N}
.$
\end{lemma}

\begin{proof}
See the proof of Proposition 4.5 in \cite{MR2112694}.
\end{proof}

\begin{theorem}
\label{Th.UniversalOmar'sLemma}Let $G=(V_{G},$ $E_{G},$ $\pi _{G})$ be a
digraph containing nontrivial strongly connected components. If the loop
number of every nontrivial strongly connected component is equal to $1$ then
there is an $m\in 
\mathbb{N}
$ such that \textbf{any} pair of vertices $a_{i},a_{j}\in V_{G}$ with $a_{i}$
recurrently connected to $a_{j}$ satisfies%
\begin{equation*}
s_{m+\lambda }(a_{i},a_{j})>0\text{ }\forall \text{ }\lambda \in 
\mathbb{N}
_{0}
\end{equation*}
\end{theorem}

\begin{proof}
Let $V_{G}=\{a_{1},...,a_{n}\}$ be the numeration of the vertices and $%
a_{i},a_{j}\in V_{G}.$ If $a_{i}$ is recurrently connected to $a_{j}$, then
necessarily there is a sequence $a_{i}\rightsquigarrow _{s}a_{j}$ that
visits a vertex contained in a nontrivial strongly connected component. In
other words, $\exists $ $a_{k}\in V_{f}$ \ and a sequence $%
a_{i}\rightsquigarrow _{s}a_{j}$ such that $\overleftrightarrow{a_{k}}$ is
nontrivial and%
\begin{equation*}
a_{i}\rightsquigarrow _{s}a_{j}=a_{i}\rightarrow ...\rightarrow
a_{k}\rightarrow ...\rightarrow a_{j}
\end{equation*}%
By Lemma \ref{Omar'sProposition} there is a $m_{k}\in 
\mathbb{N}
$ such that there are sequences  $a_{k}\rightsquigarrow
_{m_{k}+\lambda }a_{k}$ $\forall $ $\lambda \in 
\mathbb{N}
_{0}.$ Now $\forall $ $\lambda \in 
\mathbb{N}
_{0}$ we can construct a sequence%
\begin{equation*}
a_{i}\rightsquigarrow _{s_{\lambda }}a_{j}=a_{i}\rightarrow ...\rightarrow
a_{k}\rightsquigarrow _{m_{k}+\lambda }a_{k}\rightarrow ...\rightarrow a_{j}
\end{equation*}%
Now, if we consider among all pairs $i,j\in $ $\{1,...,n\}$ such that $%
a_{i}\in V_{G}$ is recurrently connected to $a_{j}\in V_{G}$ the maximum $m$
of all values $m_{k}$ we can state: $\exists $ $m\in 
\mathbb{N}
$ such that any pair of recurrently connected vertices $a_{i},a_{j}\in V_{G}$
satisfies%
\begin{equation*}
s_{m+\lambda }(a_{i},a_{j})>0\text{ }\forall \text{ }\lambda \in 
\mathbb{N}
_{0}
\end{equation*}
\end{proof}

\begin{theorem}
\label{Th.BooleanFPSequiv.LoopNumber1}Let $\mathbf{F}_{2}$ be the finite
field with two elements, $f\in MF_{n}^{n}(\mathbf{F}_{2})$ a \textbf{Boolean}%
\textit{\ }coupled monomial dynamical system and $G_{f}=(V_{f},$ $E_{f},$ $%
\pi _{f})$ its dependency graph. $f$ is a fixed point system if and only if
the loop number of each nontrivial strongly connected components of $G_{f}$
is equal to $1$.
\end{theorem}

\begin{proof}
The necessity follows from Corollary \ref{CorNecessityOfLoopNo1}. Now assume
that each nontrivial strongly connected components of $G_{f}$ has loop
number $1$ and let $V_{f}=\{a_{1},...,a_{n}\}$ be the numeration of the
vertices. Furthermore let $F:=\Psi ^{-1}(f)$ be the corresponding matrix and
consider vertices $a_{i},a_{j}\in V_{f}.$ By Remark \ref{RemarkRecurrentConn}%
, either $a_{i}$ is recurrently connected to $a_{j}$ or there is an $%
u_{0}\in 
\mathbb{N}
$ with $u_{0}\leq n$ such that no sequence $a_{i}\rightsquigarrow _{t}a_{j}$
of length $t\geq u_{0}$ exists. If the latter is the case, then%
\begin{equation*}
(F^{u_{0}+\lambda })_{ij}=0\text{ }\forall \text{ }\lambda \in 
\mathbb{N}
_{0}
\end{equation*}%
On the other hand, if $a_{i}$ is recurrently connected to $a_{j}$, then by
Theorem \ref{Th.UniversalOmar'sLemma} there is an $m_{0}\in 
\mathbb{N}
$ such that%
\begin{equation*}
(F^{m_{0}+\lambda })_{ij}\neq 0\text{ }\forall \text{ }\lambda \in 
\mathbb{N}
_{0}
\end{equation*}%
Therefore, we have for $m:=\max (m_{0},u_{0})$ that%
\begin{equation*}
(F^{m+\lambda })_{ij}\neq 0\text{ }\forall \text{ }\lambda \in 
\mathbb{N}
_{0}\text{ or }(F^{m+\lambda })_{ij}=0\text{ }\forall \text{ }\lambda \in 
\mathbb{N}
_{0}
\end{equation*}%
Summarizing we have by 2. of Lemma \ref{red.alg.Lemma}%
\begin{equation*}
mred_{q}(F^{m+\lambda })=mred_{q}(F^{m})\text{ }\forall \text{ }\lambda \in 
\mathbb{N}%
\end{equation*}%
and by \ref{RemarkPowersOfF}%
\begin{equation*}
F^{\cdot (m+\lambda )}=F^{\cdot m}\text{ }\forall \text{ }\lambda \in 
\mathbb{N}%
\end{equation*}%
Thus, after applying the isomorphism $\Psi $%
\begin{equation*}
f^{m+\lambda }=f^{m}\text{ }\forall \text{ }\lambda \in 
\mathbb{N}%
\end{equation*}
\end{proof}

\begin{remark}
The statements of the previous theorems together with the Remark \ref%
{Rem.ZeroExcludedFromMonSyst.} about zero functions as components constitute
the statement of Theorem 6.1 in \cite{MR2112694}.
\end{remark}

In the following two corollaries we provide alternative proofs to the claims
made in Corollary 6.3 and Theorem 6.5 of \cite{MR2112694}:

\begin{corollary}[and Definition]
Let $\mathbf{F}_{2}$ the finite field with two elements and $f\in MF_{n}^{n}(%
\mathbf{F}_{2})$ the coupled monomial dynamical system defined by%
\begin{eqnarray*}
f_{1}(x) &=&x_{1}^{a_{11}} \\
f_{i}(x) &=&(\prod\limits_{j=1}^{i-1}x_{j}^{a_{ij}})x_{i}^{a_{ii}},\text{ }%
i=2,...,n
\end{eqnarray*}%
where $a_{ij}\in E_{q},$ $i=1,...,n,$ $j=1,...,i-1$. Such a system is called
a \textit{Boolean triangular system}. Boolean triangular systems are always
fixed point systems.
\end{corollary}

\begin{proof}
From the structure of $f$ it is easy to see that every strongly connected
component of the dependency graph of $f$ is either trivial or has loop
number $1$.
\end{proof}

\begin{corollary}
Let $\mathbf{F}_{2}$ the finite field with two elements, $f\in MF_{n}^{n}(%
\mathbf{F}_{2})$ a fixed point system and $j,i\in \{1,...,n\}$. Consider the
system  $g\in MF_{n}^{n}(\mathbf{F}_{2})$ defined as $%
g_{k}(x):=f_{k}(x)$ $\forall $ $k\in \{1,...,n\}\backslash j$ and $%
g_{j}(x):=x_{i}f_{j}(x)$ $\forall $ $x\in \mathbf{F}_{2}^{n}.$ Then $g$ is a
fixed point system if there is no sequence $a_{i}\rightsquigarrow _{s}a_{j}$
from $a_{i}$ to $a_{j}$ or if $\overleftrightarrow{a_{i}}$ or $%
\overleftrightarrow{a_{j}}$ are nontrivial.
\end{corollary}

\begin{proof}
If $i=j$ then $E_{g}$ contains the self loop $a_{i}\rightarrow a_{i}$ and $%
\overleftrightarrow{a_{i}}$ becomes nontrivial (if it wasn't already) with
loop number $1.$ If $i\neq j$ then we have two cases: If there is no
sequence $a_{i}\rightsquigarrow _{s}a_{j},$ then adding the edge $%
a_{j}\rightarrow a_{i}$ (which might be already there) doesn't affect $%
\overleftrightarrow{a_{i}}\neq \overleftrightarrow{a_{j}}.$ If there is a
sequence $a_{i}\rightsquigarrow _{s}a_{j}$ then adding the edge $%
a_{j}\rightarrow a_{i}$ (which might be already there) forces $%
\overleftrightarrow{a_{i}}=\overleftrightarrow{a_{j}}.$ Now since by
hypothesis $\overleftrightarrow{a_{i}}$ or $\overleftrightarrow{a_{j}}$ are
nontrivial and $f$ is a fixed point system, then 
\begin{equation*}
\tciLaplace (\overleftrightarrow{a_{i}})=\tciLaplace (\overleftrightarrow{%
a_{j}})=1
\end{equation*}
\end{proof}

\begin{definition}
\label{Def.(q-1)-foldly redundant syst.}Let $\mathbf{F}_{q}$ be a finite
field, $f\in MF_{n}^{n}(\mathbf{F}_{q})$ a monomial dynamical system and $%
G_{f}=(V_{f},$ $E_{f},$ $\pi _{f})$ its dependency graph. $f$ is called a $%
(q-1)$\textit{-fold redundant monomial system} if there is an $N\in 
\mathbb{N}
$ such that for \textbf{any} pair $a,b\in $ $V_{f}$ with $a$ recurrently
connected to $b,$ the following holds:%
\begin{equation*}
\forall \text{ }m\geq N\text{ }\exists \text{ }\alpha _{abm}\in 
\mathbb{N}
_{0}:s_{m}(a,b)=\alpha _{abm}(q-1)
\end{equation*}
\end{definition}

\begin{lemma}
\label{(q-1)-foldly redundant syst.FPSLemma}Let $\mathbf{F}_{q}$ be a finite
field, $f\in MF_{n}^{n}(\mathbf{F}_{q})$ a coupled $(q-1)$-fold redundant
monomial dynamical system and $G_{f}=(V_{f},$ $E_{f},$ $\pi _{f})$ its
dependency graph. Then $f$ is a fixed point system if the loop number of
each nontrivial strongly connected component of $G_{f}$ is equal to $1.$
\end{lemma}

\begin{proof}
Let $V_{f}=\{a_{1},...,a_{n}\}$ be the numeration of the vertices and $%
F:=\Psi ^{-1}(f)$ be the corresponding matrix of $f.$ Consider two arbitrary
vertices $a_{i},a_{j}\in V_{f}.$ By Remark \ref{RemarkRecurrentConn}, either 
$a_{i}$ is recurrently connected to $a_{j}$ or there is an $m_{0}\in 
\mathbb{N}
$ with $m_{0}\leq n$ such that no sequence $a\rightsquigarrow _{t}b$ of
length $t\geq m_{0}$ exists. If the latter is the case, then%
\begin{equation*}
(F^{m_{0}+\lambda })_{ij}=0\text{ }\forall \text{ }\lambda \in 
\mathbb{N}
_{0}
\end{equation*}%
On the other hand, if $a_{i}$ is recurrently connected to $a_{j}$, then by
Theorem \ref{Th.UniversalOmar'sLemma} there is an $m_{1}\in 
\mathbb{N}
$ such that%
\begin{equation}
s_{m_{1}+\gamma }(a_{i},a_{j})>0\text{ }\forall \text{ }\gamma \in 
\mathbb{N}
_{0}  \label{eq.sPos0}
\end{equation}%
Consider now $m_{2}:=\max (n,m_{1}).$ Due to the universality of $m_{1}$ in
the expression (\ref{eq.sPos0}), for any pair of vertices $a_{i},a_{j}\in
V_{G}$ with $a_{i}$ recurrently connected to $a_{j}$ there is a sequence $%
a_{i}\rightsquigarrow _{m_{2}+\gamma }a_{j}$ of length $m_{2}+\gamma ,$ in
particular $s_{(m_{2}+\gamma )}(a_{i},a_{j})>0$ $\forall $ $\gamma \in 
\mathbb{N}
_{0}.$ Now, let $N$ be the constant in Definition \ref{Def.(q-1)-foldly
redundant syst.} and $m_{3}:=\max (N,m_{2}).$ Now, by hypothesis, $\exists $ 
$\alpha _{ij\gamma }\in 
\mathbb{N}
$ such that%
\begin{equation*}
s_{(m_{3}+\gamma )}(a_{i},a_{j})=\alpha _{ij\gamma }(q-1)\text{ }\forall 
\text{ }\gamma \in 
\mathbb{N}
_{0}
\end{equation*}%
Thus%
\begin{eqnarray*}
s_{(m_{3}+\gamma )}(a_{i},a_{j}) &=&\alpha _{ij\gamma }(q-1)=\alpha
_{ij0}(q-1)+(\alpha _{ij\gamma }-\alpha _{ij0})(q-1) \\
&=&s_{m_{3}}(a_{i},a_{j})+(\alpha _{ij\gamma }-\alpha _{ij0})(q-1)\text{ }%
\forall \text{ }\gamma \in 
\mathbb{N}
_{0}
\end{eqnarray*}%
Summarizing, since $m_{0}\leq n\leq m_{2}\leq m_{3},$ we can say $\forall $ $%
i,j\in $ $\{1,...,n\}$, depending on whether $a_{i}$ and $a_{j}$ are
recurrently connected or not,%
\begin{equation*}
(F^{m_{3}+\lambda })_{ij}=0\text{ }\forall \text{ }\lambda \in 
\mathbb{N}
_{0}
\end{equation*}%
or%
\begin{equation*}
\exists \text{ }a_{\lambda }\in 
\mathbb{Z}
:(F^{m_{3}+\lambda })_{ij}=(F^{m_{3}})_{ij}+a_{\lambda }(q-1)\neq 0\text{ }%
\forall \text{ }\lambda \in 
\mathbb{N}
_{0}
\end{equation*}%
Now, by 2. and 4. of Lemma \ref{red.alg.Lemma} it follows%
\begin{equation*}
mred_{q}(F^{m_{3}+\lambda })=mred_{q}(F^{m_{3}})\text{ }\forall \text{ }%
\lambda \in 
\mathbb{N}%
\end{equation*}%
and by \ref{RemarkPowersOfF}%
\begin{equation*}
F^{\cdot (m_{3}+\lambda )}=F^{\cdot m_{3}}\text{ }\forall \text{ }\lambda
\in 
\mathbb{N}%
\end{equation*}%
Thus, after applying the isomorphism $\Psi $%
\begin{equation*}
f^{m_{3}+\lambda }=f^{m_{3}}\text{ }\forall \text{ }\lambda \in 
\mathbb{N}%
\end{equation*}
\end{proof}

\begin{theorem}
Let $\mathbf{F}_{q}$ be a finite field, $f\in MF_{n}^{n}(\mathbf{F}_{q})$ a
coupled monomial dynamical system and $G_{f}=(V_{f},$ $E_{f},$ $\pi _{f})$
its dependency graph. Then $f$ is a fixed point system if the following
properties hold

\begin{enumerate}
\item {\normalsize The loop number of each nontrivial strongly connected
component of }$G_{f}${\normalsize \ is equal to }$1.$

\item {\normalsize For each nontrivial strongly connected component }$%
\overleftrightarrow{a}\subseteq V_{f}${\normalsize \ and arbitrary }$b,c\in $%
{\normalsize \ }$\overleftrightarrow{a}${\normalsize ,}%
\begin{equation*}
s_{1}(b,c)\neq 0\Rightarrow s_{1}(b,c)=q-1
\end{equation*}
\end{enumerate}
\end{theorem}

\begin{proof}
Let $V_{f}=\{a_{1},...,a_{n}\}$ be the numeration of the vertices and $%
F:=\Psi ^{-1}(f)$ be the corresponding matrix of $f.$ Consider two
vertices  $a_{i},a_{j}\in V_{f}$ such that $a_{i}$ is recurrently
connected to $a_{j}$. Then by Theorem \ref{Th.UniversalOmar'sLemma} there is
an $m_{1}\in 
\mathbb{N}
$ such that%
\begin{equation}
s_{m_{1}+\gamma }(a_{i},a_{j})>0\text{ }\forall \text{ }\gamma \in 
\mathbb{N}
_{0}  \label{eq.sPos}
\end{equation}%
Consider now $m_{2}:=\max (n,m_{1}).$ Due to the universality of $m_{1}$ in
the expression (\ref{eq.sPos}), for any pair of vertices $a_{i},a_{j}\in
V_{G}$ with $a_{i}$ recurrently connected to $a_{j}$ there is a sequence $%
a_{i}\rightsquigarrow _{m_{2}+\gamma }a_{j}$ of length $m_{2}+\gamma .$
Since $m_{2}+\gamma >n-1,$ necessarily $\exists $ $a_{k_{\gamma
}},a_{l_{\gamma }}\in \overleftrightarrow{a_{k_{\gamma }}}$ such that $%
\overleftrightarrow{a_{k_{\gamma }}}$ is nontrivial and%
\begin{equation}
a_{i}\rightsquigarrow _{(m_{2}+\gamma )}a_{j}=a_{i}\rightarrow
...\rightarrow a_{k_{\gamma }}\rightsquigarrow _{t}a_{l_{\gamma
}}\rightarrow ...\rightarrow a_{j}  \label{eq.ssequence}
\end{equation}%
($t$ depends on $i,j$ and $\gamma $). Now, by hypothesis, every two directly
connected vertices $a,b\in \overleftrightarrow{a_{k_{\gamma }}}$ are
directly connected by exactly $q-1$ directed edges. Therefore, for any
sequence $a_{k_{\gamma }}\rightsquigarrow _{t}a_{l_{\gamma }}$ of length $%
t\in 
\mathbb{N}
$ there are $(q-1)^{t}$ different copies of it and we can conclude$\exists $ 
$\alpha \in 
\mathbb{N}
$ such that $s_{t}(a_{k_{\gamma }},a_{l_{\gamma }})=\alpha (q-1).$ As a
consequence, there are $\alpha (q-1)$ different copies of the sequence (\ref%
{eq.ssequence}). Since we are dealing with an arbitrary sequence $%
a_{i}\rightsquigarrow _{(m_{2}+\gamma )}a_{j}$ of fixed length $m_{2}+\gamma
,$ $\gamma \in 
\mathbb{N}
_{0}$ we can conclude that $\exists $ $\alpha _{ij\gamma }\in 
\mathbb{N}
$ such that%
\begin{equation*}
s_{(m_{2}+\gamma )}(a_{i},a_{j})=\alpha _{ij\gamma }(q-1)\text{ }\forall 
\text{ }\gamma \in 
\mathbb{N}
_{0}
\end{equation*}%
Thus $f$ is a coupled $(q-1)$-fold redundant monomial dynamical system and
the claim follows from Lemma \ref{(q-1)-foldly redundant syst.FPSLemma}.
\end{proof}

\begin{corollary}
Let $\mathbf{F}_{2}$ be the finite field with two elements, $f\in MF_{n}^{n}(%
\mathbf{F}_{2})$ a \textit{Boolean} monomial dynamical system and $F:=\Psi
^{-1}(f)\in M(n\times n;$ $E_{2})$ its corresponding matrix. Furthermore,
let $\mathbf{F}_{q}$ be a finite field and $g\in MF_{n}^{n}(\mathbf{F}_{q})$
the monomial dynamical system whose corresponding matrix $G:=\Psi
^{-1}(g)\in M(n\times n;$ $E_{q})$ satisfies $\forall $ $i,j\in \{1,...,n\}$%
\begin{equation*}
G_{ij}=\left\{ 
\begin{array}{c}
q-1\text{ if }F_{ij}=1 \\ 
0\text{ if }F_{ij}=0%
\end{array}%
\right.
\end{equation*}%
If $f$ is a fixed point system then $g$ is a fixed point system too.
\end{corollary}

\begin{proof}
Let $G_{f}=(V_{f},$ $E_{f},$ $\pi _{f})$ be the dependency graph of $f.$ By
the definition of $g,$ one can easily see that the dependency
graph  $G_{g}=(V_{g},$ $E_{g},$ $\pi _{g})$ of $g$ can be generated
from $G_{f}$ by adding $q-2$ identical parallel edges for every existing
edge. Obviously $G_{f}$ and $G_{g}$ have the same strongly connected
components. If $G_{f}$ doesn't contain any nontrivial strongly connected
components, then $G_{g}$ wouldn't contain any either and by Theorem \ref%
{Th.FPSwhenOnlyTrivialSCC} $g$ would be a fixed point system. If, on the
other hand, $G_{f}$ does contain nontrivial strongly connected components,
then by Theorem \ref{Th.BooleanFPSequiv.LoopNumber1} each of those
components would have loop number $1.$ From the definition of $g$ it also
follows for any pair of vertices $a,b\in E_{g}$%
\begin{equation*}
s_{1}(a,b)\neq 0\Rightarrow s_{1}(a,b)=q-1
\end{equation*}%
By the previous theorem $g$ would be a fixed point system.
\end{proof}

\begin{example}[and Corollary]
Let $\mathbf{F}_{q}$ be a finite field and $f\in MF_{n}^{n}(\mathbf{F}_{q})$
the coupled monomial dynamical system defined by%
\begin{eqnarray*}
f_{1}(x) &=&x_{1}^{q-1} \\
f_{i}(x) &=&(\prod\limits_{j=1}^{i-1}x_{j}^{a_{ij}})x_{i}^{q-1},\text{ }%
i=2,...,n
\end{eqnarray*}%
where $a_{ij}\in E_{q},$ $i=1,...,n,$ $j=1,...,i-1$ are not further
specified exponents. Such a system is called \textit{triangular}. It is easy
to see that the dependency graph of $f$ contains $n$ one vertex nontrivial
strongly connected components. Each of them has a $(q-1)$-fold self loop.
Therefore, by the previous Theorem, $f$ must be a fixed point system.
\end{example}

\begin{theorem}
Let $\mathbf{F}_{q}$ be a finite field, $f\in MF_{n}^{n}(\mathbf{F}_{q})$ a
coupled monomial dynamical system and $G_{f}=(V_{f},$ $E_{f},$ $\pi _{f})$
its dependency graph. Then $f$ is a fixed point system if for every vertex $%
a\in V_{f}$ that is recurrently connected to some other vertex $b\in V_{f}$
the edge $a\rightarrow a$ appears exactly $q-1$ times in $E_{f},$ i.e.%
\begin{equation*}
\left\vert \pi _{f}^{-1}((a,a))\right\vert =q-1
\end{equation*}
\end{theorem}

\begin{proof}
Let $V_{f}=\{a_{1},...,a_{n}\}$ be the numeration of the vertices and $%
F:=\Psi ^{-1}(f)$ be the corresponding matrix of $f.$ Consider two
vertices  $a_{i},a_{j}\in V_{f}$ such that $a_{i}$ is recurrently
connected to $a_{j}$. Then by Theorem \ref{Th.UniversalOmar'sLemma} there is
an $m_{1}\in 
\mathbb{N}
$ such that%
\begin{equation}
s_{m_{1}+\gamma }(a_{i},a_{j})>0\text{ }\forall \text{ }\gamma \in 
\mathbb{N}
_{0}  \label{eq.sPos2}
\end{equation}%
Consider now $m_{2}:=\max (n,m_{1}).$ Due to the universality of $m_{1}$ in
the expression (\ref{eq.sPos2}), for any pair of vertices $a_{i},a_{j}\in
V_{G}$ with $a_{i}$ recurrently connected to $a_{j}$ there is a sequence $%
a_{i}\rightsquigarrow _{m_{2}+\gamma }a_{j}$ of length $m_{2}+\gamma .$
Consider one particular sequence $a_{i}\rightsquigarrow _{m_{2}+\gamma
}a_{j} $ of length $m_{2}+\gamma $ and call it $w_{\gamma
}:=a_{i}\rightsquigarrow _{m_{2}+\gamma }a_{j}.$ By hypothesis there are
exactly $q-1$ directed edges $a_{i}\rightarrow a_{i}$. Therefore, there are $%
q-1$ copies of the sequence $w_{\gamma }.$ Since we are dealing with an
arbitrary sequence $a_{i}\rightsquigarrow _{(m_{2}+\gamma )}a_{j}$ of fixed
length $m_{2}+\gamma ,$ $\gamma \in 
\mathbb{N}
_{0}$ we can conclude that $\exists $ $\alpha _{ij\gamma }\in 
\mathbb{N}
$ such that%
\begin{equation*}
s_{(m_{2}+\gamma )}(a_{i},a_{j})=\alpha _{ij\gamma }(q-1)\text{ }\forall 
\text{ }\gamma \in 
\mathbb{N}
_{0}
\end{equation*}%
Thus $f$ is a coupled $(q-1)$-fold redundant monomial dynamical system and
the claim follows from Lemma \ref{(q-1)-foldly redundant syst.FPSLemma}.
\end{proof}

\begin{example}[and Corollary]
Let $\mathbf{F}_{q}$ be a finite field and $f\in MF_{n}^{n}(\mathbf{F}_{q})$
a monomial dynamical system such that the diagonal entries of its
corresponding matrix $F:=\Psi ^{-1}(f)$ satisfy%
\begin{equation*}
F_{ii}=q-1\text{ }\forall \text{ }i\in \{1,...,n\}
\end{equation*}%
Since every vertex satisfies the requirement of the previous theorem, $f$
must be a fixed point system. This result generalizes our previous result
about triangular monomial dynamical systems.
\end{example}

We now provide a more algebraic sufficient condition for a system $f\in
MF_{n}^{n}(\mathbf{F}_{q})$ to be a fixed point system.

\begin{lemma}
Let $n\in 
\mathbb{N}
$ be a natural number and $A\in M(n\times n;$ $%
\mathbb{R}
)$ a real matrix. In addition, let $A$ be diagonalizable over $%
\mathbb{C}
.$ Then $A^{m}=A$ $\forall $ $m\in 
\mathbb{N}
$ if and only if $\exists $ $r,s\in 
\mathbb{N}
_{0}$ such that $r+s=n$ and the characteristic polynomial $charpoly(A)$ of $%
A $ can be written as%
\begin{equation*}
charpoly(A)=a(\lambda -1)^{s}\lambda ^{t}
\end{equation*}%
where $a\in 
\mathbb{R}
\backslash \{0\}.$
\end{lemma}

\begin{proof}
The proof of this simple linear algebraic result is left to the interested
reader.
\end{proof}

\begin{theorem}
Let $\mathbf{F}_{q}$ be a finite field, $f\in MF_{n}^{n}(\mathbf{F}_{q})$ a
coupled monomial dynamical system and $F:=\Psi ^{-1}(f)\in M(n\times n;$ $%
E_{q})$ its corresponding matrix. If the matrix $F$ (viewed as a real matrix 
$F$ $\in M(n\times n;$ $%
\mathbb{N}
)\subset M(n\times n;$ $%
\mathbb{R}
)$) has the characteristic polynomial%
\begin{equation}
charpoly(F)=a(\lambda -1)^{s}\lambda ^{t}  \label{eq.CharPolyCond.}
\end{equation}%
where $a\in 
\mathbb{Z}
\backslash \{0\},$ $r,s\in 
\mathbb{N}
_{0}$ such that $r+s=n$ and the geometric multiplicity of the eigenvalues $0$
and $1$ is equal to the corresponding algebraic multiplicity, then $f$ is a
fixed point system.
\end{theorem}

\begin{proof}
It is a well-known linear algebraic result that if there is a basis of
eigenvectors of a matrix, the matrix is diagonalizable. By the hypothesis
this is the case for $F$. Therefore, by the previous Lemma%
\begin{equation*}
F^{m}=F\text{ }\forall \text{ }m\in 
\mathbb{N}%
\end{equation*}%
Now, by Remarks \ref{RemarkMatPower} and \ref{RemarkMatRed} we consequently
have $\forall $ $m\in 
\mathbb{N}
$%
\begin{equation*}
\Psi ^{-1}(f^{m})=F^{\cdot m}=mred_{q}(F^{m})=mred_{q}(F)=F
\end{equation*}%
After applying the isomorphism $\Psi $ we get%
\begin{equation*}
f^{m}=f\text{ \ }\forall \text{ }m\in 
\mathbb{N}%
\end{equation*}
\end{proof}

\begin{remark}
Let $\mathbf{F}_{q}$ be a finite field, $f\in MF_{n}^{n}(\mathbf{F}_{q})$ a
coupled monomial dynamical system and $F:=\Psi ^{-1}(f)\in M(n\times n;$ $%
E_{q})$ its corresponding matrix. The matrix $F$ viewed as the adjacency
matrix of the dependency graph $G_{f}=(V_{f},$ $E_{f},$ $\pi _{f})$ of $f$
satisfies%
\begin{equation*}
F^{m}=F\text{ }\forall \text{ }m\in 
\mathbb{N}%
\end{equation*}%
if and only if for each pair of vertices $a,b\in V_{f}$ \ the value $%
s_{m}(a,b)$ is constant for all $m\in 
\mathbb{N}
$. In other words, $a$ and $b$ are either disconnected or for every length $%
m\in 
\mathbb{N}
$ they are connected with the same degree of redundancy.
\end{remark}

\begin{example}
Consider the monomial system $g\in MF_{5}^{5}(\mathbf{F}_{3})$ defined by
the matrix%
\begin{equation*}
G:=%
\begin{pmatrix}
1 & 1 & 0 & 0 & 0 \\ 
0 & 1 & 1 & 0 & 0 \\ 
0 & 0 & 1 & 0 & 0 \\ 
0 & 0 & 0 & 0 & 1 \\ 
0 & 0 & 0 & 0 & 0%
\end{pmatrix}%
\end{equation*}%
It is easy to show that%
\begin{equation*}
charpoly(G)=(\lambda -1)^{3}\lambda ^{2}
\end{equation*}%
However, $g$ is not a fixed point system. This shows that the condition (\ref%
{eq.CharPolyCond.}) alone is not sufficient.
\end{example}

\section{An algorithm of polynomial complexity to identify fixed point
systems}

\subsection{Some basic considerations}

\begin{definition}
Let $X$ be a nonempty finite set, $n\in 
\mathbb{N}
$ a natural number and $f:X^{n}\rightarrow X^{n}$ a time discrete finite
dynamical system. The \textit{phase space of }$f$ is the digraph with node
set $X^{n}$, arrow set $E$ defined as%
\begin{equation*}
E:=\{(x,y)\in X^{n}\times X^{n}\text{ }|\text{ }f(x)=y\}
\end{equation*}%
and vertex mapping%
\begin{eqnarray*}
\pi &:&E\rightarrow X^{n}\times X^{n} \\
(x,y) &\mapsto &(x,y)
\end{eqnarray*}
\end{definition}

\begin{remark}
Due to the finiteness of $X$ it is obvious that the trajectory%
\begin{equation*}
x,\text{ }f(x),\text{ }f^{2}(x),...
\end{equation*}%
of any point $x\in X^{n}$ contains at most $\left\vert X^{n}\right\vert
=\left\vert X\right\vert ^{n}$ different points and therefore becomes either
cyclic or converges to a single point $y\in X$ with the property $f(y)=y$
(i.e. a fixed point of $f$). Thus, the phase space consists of closed paths
of different lengths between $1$ (i.e. loops centered on fixed points) and $%
\left\vert X^{n}\right\vert =\left\vert X\right\vert ^{n}$ and directed
trees that end each one at exactly one closed path. The nodes in the
directed trees correspond to \textit{transient states} of the system.
\end{remark}

According to our definition of monomial dynamical system $f\in MF_{n}^{n}(%
\mathbf{F}_{q})$, the possibility that one of the functions $f_{i}$ is equal
to the zero function is excluded (see Definition \ref{Def.MonDynSyst} and
Remark \ref{Rem.ZeroExcludedFromMonSyst.}). Therefore, the following
algorithm is designed for such systems. However, in this algorithmic
framework it would be convenient to include the more general case (as
defined in \cite{MR2293353} and \cite{MR2112694}), i.e. the case when some
of the functions $f_{i}$ can indeed be equal to the zero function. In the
vein of Remark \ref{Rem.ZeroExcludedFromMonSyst.} this actually only
requires some type of preprocessing. The preprocessing algorithm will be
described and analyzed in the Appendix.

Our algorithm is based on the following observation made by Dr. Michael
Shapiro about general time discrete finite dynamical systems: By the
previous remark, a chain of transient states in the phase space of a time
discrete finite dynamical system $f:X^{n}\rightarrow X^{n}$ can contain at
most $s:=\left\vert X^{n}\right\vert -1=\left\vert X\right\vert ^{n}-1$
transient elements. Therefore, to determine whether a system is a fixed
point system, it is sufficient to establish whether the mappings $f^{r}$ and 
$f^{r+1}$ are identical for any $r\geq s.$ In the case of a monomial system $%
f\in MF_{n}^{n}(\mathbf{F}_{q}),$ due to Theorem \ref{ThmMonoidIsomorphism},
we only need to look at the corresponding matrices  $F^{\cdot r},$ $%
F^{\cdot r+1}\in M(n\times n;E_{q}).$ Computationally it is more convenient
to generate the following sequence of powers%
\begin{equation*}
F^{\cdot 2},\text{ }(F^{\cdot 2})^{\cdot 2}=F^{\cdot 4},\text{ }(F^{\cdot
4})^{\cdot 2}=F^{\cdot 8},\text{ }(F^{\cdot 8})^{\cdot 2}=F^{\cdot
16},...,F^{\cdot (2^{t})}
\end{equation*}%
To achieve the "safe" number of iterations $\left\vert \mathbf{F}%
_{q}^{n}\right\vert -1=q^{n}-1$ we need to make sure%
\begin{equation*}
2^{t}\geq q^{n}-1
\end{equation*}%
This is equivalent to%
\begin{equation*}
t\geq \log _{2}(q^{n}-1)
\end{equation*}%
To obtain a natural number we use the ceil function%
\begin{equation}
t:=ceil(\log _{2}(q^{n}-1))  \label{eq.Def.t}
\end{equation}%
Thus we have, due to the monotonicity of the $\log $ function,%
\begin{equation*}
t<\log _{2}(q^{n}-1)+1\leq \log _{2}(q^{n})+1=n\log _{2}(q)+1
\end{equation*}

\subsection{The algorithm and its complexity analysis}

The algorithm is fairly simple: Given a monomial system $f\in MF_{n}^{n}(%
\mathbf{F}_{q})$ and its corresponding matrix $F:=\Psi ^{-1}(f)\in M(n\times
n;E_{q})$

\begin{enumerate}
\item {\normalsize With }$t${\normalsize \ as defined above (\ref{eq.Def.t}%
), calculate the matrices }$A:=${\normalsize \ }$F^{\cdot 2^{t}}$%
{\normalsize \ and }$B:=FA.${\normalsize \ This step requires }$t+1$%
{\normalsize \ matrix multiplications.}

\item {\normalsize Compare the }$n^{2}${\normalsize \ entries }$A_{ij}$%
{\normalsize \ and }$B_{ij}.${\normalsize \ This step requires at most }$%
n^{2}${\normalsize \ comparisons. (This maximal value is needed in the case
that }$f${\normalsize \ is a fixed point system).}

\item $f${\normalsize \ is a fixed point system if and only if the matrices }%
$A${\normalsize \ and }$B${\normalsize \ are equal.}
\end{enumerate}

It is well known that matrix multiplication requires $2n^{3}-n^{2}$ addition
or multiplication operations. Since $t+1<n\log _{2}(q)+2,$ the number of
operations required in step $1$ is bounded above by%
\begin{equation*}
(2n^{3}-n^{2})(n\log _{2}(q)+2)
\end{equation*}%
Summarizing, we have the following upper bound $N(n,q)$ for the number of
operations in steps $1$ and $2$%
\begin{equation*}
N(n,q):=(2n^{3}-n^{2})(n\log _{2}(q)+2)+n^{2}
\end{equation*}%
For a fixed size $q$ of the finite field $\mathbf{F}_{q}$ used it holds%
\begin{equation*}
\lim_{n\rightarrow \infty }\frac{N(n,q)}{n^{4}}=2\log _{2}(q)
\end{equation*}%
and we can conclude $N(n,q)\in O($ $n^{4})$ for a fixed $q.$ The asymptotic
behavior for a growing number of variables and growing number of field
elements is described by%
\begin{equation*}
\lim_{\substack{ n\rightarrow \infty  \\ q\rightarrow \infty }}\frac{N(n,q)}{%
n^{4}\log _{2}(q)}=2
\end{equation*}%
Thus, $N(n,q)\in O($ $n^{4}\log _{2}(q))$ for $n,q\rightarrow \infty .$%
\newline
It is pertinent to comment on the arithmetic operations performed during the
matrix multiplications. Since the matrices are elements of the matrix monoid 
$M(n\times n;E_{q}),$ the arithmetic operations are operations in the monoid 
$E_{q}.$ By the Lemmas \ref{AddLemma} and \ref{MultLemma} the addition resp.
the multiplication operation on $E_{q}$ requires an integer number addition%
\footnote{%
See Chapter 4 of \cite{Kaplan} for a detailed description of integer number
representation and arithmetic in typical computer algebra systems.} resp.
multiplication and a reduction as defined in Lemma \ref{red.alg.Lemma}. The
reduction $red_{q}(a)$ of an integer number $a\in 
\mathbb{N}
_{0},$ $a\geq q$ is obtained as the degree of the remainder of the
polynomial division $\tau ^{a}\div ($ $\tau ^{q}-\tau ).$ According to 4.6.5
of \cite{Kaplan} this division requires%
\begin{equation*}
O(2(\deg (\tau ^{a})-\deg (\tau ^{q}-\tau )))=O(2(a-q))
\end{equation*}%
integer number operations. However, we know that the reductions $red_{q}(.)$
are applied to the result of (regular integer) addition or multiplication of
elements of $E_{q}$ and therefore%
\begin{equation*}
a-q\leq \left\{ 
\begin{array}{c}
2(q-1)-q=q-2 \\ 
(q-1)^{2}-q=q^{2}-q+1%
\end{array}%
\right.
\end{equation*}%
As a consequence, in the worst case scenario, one addition resp.
multiplication in the monoid $E_{q}$ requires $O(q)$ resp. $O(q^{2})$
regular integer number operations.\newline
Since $E_{q}$ is a finite set and only the results of $n^{2}$ pairwise
additions and $n^{2}$ pairwise multiplications are needed, while the
algorithm is running, these numbers are of course stored in a table after
the first time they are calculated.

\section*{Acknowledgments}

We would like to thank Dr.\ Omar Col\'{o}n-Reyes for his hospitality and a
very fruitful academic interaction at the University of Puerto Rico, Mayag%
\"{u}ez. We are grateful to Prof.\ Dr.\ Bodo Pareigis for offering the
opportunity to participate in a great seminar at the
Ludwig-Maximilians-Universit\"{a}t in Munich, Germany. We also would like to
express our gratitude to Dr.\ Michael Shapiro at Tufts University, Boston,
for very helpful comments and to Dr.\ Karen Duca and Dr.\ David
Thorley-Lawson for their support. Moreover, we thank Jill Roughan, Dr.\
Michael Shapiro and Dr.\ David Thorley-Lawson for proofreading the
manuscript.

The author acknowledges support by a Public Health Service grant (RO1
AI062989) to Dr.\ David Thorley-Lawson at Tufts University, Boston, MA, USA.

\appendix

\section{Appendix}

\subsection{Some simple algebraic results}

\begin{lemma}
\label{MultLemma}Let $\mathbf{F}_{q}$ be a finite field and $a,b\in 
\mathbb{N}
_{0}$ nonnegative integers. Then it holds%
\begin{equation*}
red_{q}(ab)=red_{q}(red_{q}(a)red_{q}(b))
\end{equation*}
\end{lemma}

\begin{proof}
We have $\forall $ $x\in \mathbf{F}_{q}$%
\begin{equation*}
x^{ab}=(x^{a})^{b}=(x^{red_{q}(a)})^{red_{q}(b)}=x^{red_{q}(a)red_{q}(b)}
\end{equation*}%
and by Lemma \ref{red.alg.Lemma}%
\begin{equation*}
red_{q}(ab)=red_{q}(red_{q}(a)red_{q}(b))
\end{equation*}
\end{proof}

\begin{lemma}
\label{AddLemma}Let $\mathbf{F}_{q}$ be a finite field and $a,b\in 
\mathbb{N}
_{0}$ nonnegative integers. Then it holds%
\begin{equation*}
red_{q}(a+b)=red_{q}(red_{q}(a)+red_{q}(b))
\end{equation*}
\end{lemma}

\begin{proof}
By the division algorithm $\exists _{1}$ $g_{a}$ $,\ g_{b}$ $,$ $g_{a+b}$ $,$
$r_{a}$ $,\ r_{b}$ $,$ $r_{a+b}\in \mathbf{F}_{q}[\tau ]$ such that%
\begin{eqnarray*}
\tau ^{a} &=&g_{a}(\tau ^{q}-\tau )+r_{a}=g_{a}(\tau ^{q}-\tau )+\tau
^{red_{q}(a)} \\
\tau ^{b} &=&g_{b}(\tau ^{q}-\tau )+r_{b}=g_{b}(\tau ^{q}-\tau )+\tau
^{red_{q}(b)} \\
\tau ^{a+b} &=&g_{a+b}(\tau ^{q}-\tau )+r_{a+b}=g_{a+b}(\tau ^{q}-\tau
)+\tau ^{red_{q}(a+b)}
\end{eqnarray*}%
From the first two equations follows%
\begin{equation*}
\tau ^{a+b}=g_{a}g_{b}(\tau ^{q}-\tau )^{2}+g_{a}r_{b}(\tau ^{q}-\tau
)+r_{a}g_{b}(\tau ^{q}-\tau )+\tau ^{red_{q}(a)+red_{q}(b)}
\end{equation*}%
Applying the division algorithm to $\tau ^{red_{q}(a)+red_{q}(b)}$ we can
say  $\exists _{1}$ $g_{r}$ $,$ $r_{r}$ $\in \mathbf{F}_{q}[\tau ]$
such that%
\begin{eqnarray*}
\tau ^{a+b} &=&g_{a}g_{b}(\tau ^{q}-\tau )^{2}+g_{a}r_{b}(\tau ^{q}-\tau
)+r_{a}g_{b}(\tau ^{q}-\tau )+g_{r}(\tau ^{q}-\tau )+r_{r} \\
&=&(g_{a}g_{b}(\tau ^{q}-\tau )+g_{a}r_{b}+r_{a}g_{b}+g_{r})(\tau ^{q}-\tau
)+\tau ^{red_{q}(red_{q}(a)+red_{q}(b))}
\end{eqnarray*}%
From the uniqueness of quotient and remainder it follows%
\begin{equation*}
\tau ^{red_{q}(a+b)}=\tau ^{red_{q}(red_{q}(a)+red_{q}(b))}
\end{equation*}%
and consequently 
\begin{equation*}
red_{q}(a+b)=red_{q}(red_{q}(a)+red_{q}(b))
\end{equation*}
\end{proof}

\begin{theorem}[and Definition]
\label{Th.ExponentsSemiring}Let $\mathbf{F}_{q}$ be a finite field. The set%
\begin{equation*}
E_{q}=\{0,1,...,(q-2),(q-1)\}\subset 
\mathbb{Z}%
\end{equation*}%
together with the operations of addition $a\oplus b:=red_{q}(a+b)$ and
multiplication $a\bullet b:=red_{q}(ab)$ is a commutative semiring with
identity $1$. We call this commutative semiring the exponents semiring of
the field $\mathbf{F}_{q}.$
\end{theorem}

\begin{proof}
First we show that $E_{q}$ is a commutative monoid with respect to the
addition $\oplus $. The reduction modulo the ideal $\left\langle \tau
^{q}-\tau \right\rangle $ ensures that $E_{q}$ is closed under this
operation. Additive commutativity follows trivially from the definition. The
associativity is easily shown using Lemma \ref{AddLemma} and the fact that $%
c\in E_{q}\Leftrightarrow c=red_{q}(c)$. It is trivial to see that $0$ is
the additive identity element. $E_{q}$ is also a commutative monoid with
respect to the multiplication $\bullet :$ The reduction modulo the ideal $%
\left\langle \tau ^{q}-\tau \right\rangle $ ensures that $E_{q}$ is closed
under this operation. Multiplicative commutativity as well as the fact that $%
1$ is the multiplicative identity follow trivially from the definition. The
associativity is shown using Lemma \ref{MultLemma} and the fact that $c\in
E_{q}\Leftrightarrow c=red_{q}(c)$. The proof of the distributivity is a
straightforward verification.
\end{proof}

\begin{lemma}
Let $n\in 
\mathbb{N}
$ be a natural number, $\mathbf{F}_{q}$ be a finite field and $E_{q}$ the
exponents semiring of $\mathbf{F}_{q}.$ The set $M(n\times n;$ $E_{q})$ of $%
n\times n$ quadratic matrices with entries in the semiring $E_{q}$ together
with the operation $\cdot $ of matrix multiplication over $E_{q}$ is a
monoid.
\end{lemma}

\begin{proof}
The matrix multiplication $\cdot $ is defined in terms of the operations $%
\oplus $ and $\bullet $ on the matrix entries, therefore $M(n\times n;$ $%
E_{q})$ is closed under multiplication. The proof of the associativity is a
tedious but straightforward verification. The identity element is obviously
the unit matrix $I$.
\end{proof}

\begin{remark}
\label{RemarkMatRed}Since the entries for the matrix product $D=A\cdot B$
are defined as%
\begin{equation*}
D_{ij}=A_{i1}\bullet B_{1j}\oplus A_{i2}\bullet B_{2j}\oplus ...\oplus
A_{in}\bullet B_{nj}
\end{equation*}%
according to the definitions of the operations $\bullet $ and $\oplus $ we
can write%
\begin{eqnarray*}
D_{ij} &=&red_{q}(A_{i1}B_{1j})\oplus red_{q}(A_{i2}B_{2j})\oplus ...\oplus
red_{q}(A_{in}B_{nj}) \\
&=&red_{q}(red_{q}(A_{i1}B_{1j})+red_{q}(A_{i2}B_{2j})+...+red_{q}(A_{in}B_{nj}))
\end{eqnarray*}%
Now, by Lemma \ref{AddLemma} we have%
\begin{equation*}
D_{ij}=red_{q}(A_{i1}B_{1j}+A_{i2}B_{2j}+...+A_{in}B_{nj})
\end{equation*}%
As a consequence, if we define the following reduction operation for
matrices with nonnegative integer entries%
\begin{eqnarray*}
mred_{q} &:&M(n\times n;%
\mathbb{N}
_{0})\rightarrow M(n\times n;E_{q}) \\
A_{ij} &\mapsto &red_{q}(A_{ij})
\end{eqnarray*}%
then the following property holds for $U,V\in M(n\times n;%
\mathbb{N}
_{0})$ and  $W:=UV\in M(n\times n;%
\mathbb{N}
_{0})$%
\begin{equation*}
mred_{q}(UV)=mred_{q}(U)\cdot mred_{q}(V)
\end{equation*}%
It can be easily shown that $M(n\times n;%
\mathbb{N}
_{0})$ is a monoid and $mred_{q}:M(n\times n;%
\mathbb{N}
_{0})\rightarrow M(n\times n;E_{q})$ a monoid homomorphism. In addition, by
2. of Lemma \ref{red.alg.Lemma} we can conclude%
\begin{equation}
mred_{q}(A)=0\Leftrightarrow A=0  \label{eq.mred(A)=0}
\end{equation}
\end{remark}

\subsection{The preprocessing algorithm}

We start with the definition of monomial dynamical system according to \cite%
{MR2112694} and \cite{MR2293353}:

\begin{definition}
Let $\mathbf{F}_{q}$ be a finite field. A map $f:\mathbf{F}%
_{q}^{n}\rightarrow \mathbf{F}_{q}^{n}$ is called a monomial dynamical
system over $\mathbf{F}_{q}$ if for every $i\in \{1,...,n\}$ there exists a
tuple $(F_{i1},...,F_{in})\in E_{q}^{n}$ and an element $a_{i}\in
\{0,1\}\subseteq \mathbf{F}_{q}$ such that%
\begin{equation*}
f_{i}(x)=a_{i}x_{1}^{F_{i1}}...x_{n}^{F_{in}}\text{ }\forall \text{ }x\in 
\mathbf{F}_{q}^{n}
\end{equation*}
\end{definition}

In order to use the algorithm described in Section 4 to determine whether
such a monomial dynamical system is a fixed point system we need to
preprocess the system in the sense of Remark \ref%
{Rem.ZeroExcludedFromMonSyst.}. To accomplish this task algorithmically, we
add an element $-\infty $ to our exponents semiring $E_{q.}$ (See Definition %
\ref{Def.ExpSet} and Theorem \ref{Th.ExponentsSemiring}.):%
\begin{equation*}
\overline{E_{q}}:=E_{q}\cup \{-\infty \}
\end{equation*}%
The arithmetic with this new element is as follows%
\begin{eqnarray*}
a\oplus -\infty &=&-\infty \oplus a=-\infty \text{ }\forall \text{ }a\in 
\overline{E_{q}} \\
a\bullet -\infty &=&-\infty \bullet a=-\infty \text{ }\forall \text{ }a\in 
\overline{E_{q}}\backslash \{0\} \\
0\bullet a &=&a\bullet 0=0\text{ }\forall \text{ }a\in \overline{E_{q}}
\end{eqnarray*}%
The addition is due to the additive "absorption property" of $-\infty $
obviously associative. The same holds for the multiplication, since both $0$
and $-\infty $ show the multiplicative "absorption property" (although $0$
wins over $-\infty $). With this rules we are already able to multiply pairs
of matrices with entries in $\overline{E_{q}}.$ With this extended exponents
set we can represent the monomial dynamical systems defined above as follows:

\begin{definition}
Let $\mathbf{F}_{q}$ be a finite field. A map $f:\mathbf{F}%
_{q}^{n}\rightarrow \mathbf{F}_{q}^{n}$ is called a monomial dynamical
system over $\mathbf{F}_{q}$ if for every $i\in \{1,...,n\}$ there exists a
tuple $(F_{i1},...,F_{in})\in E_{q}^{n}$ or a tupel $(F_{i1},...,F_{in})\in
\{-\infty \}^{n}$ such that%
\begin{equation*}
f_{i}(x)=x_{1}^{F_{i1}}...x_{n}^{F_{in}}\text{ }\forall \text{ }x\in \mathbf{%
F}_{q}^{n}
\end{equation*}%
Now we describe the preprocessing algorithm: Given a monomial system $f\in
MF_{n}^{n}(\mathbf{F}_{q})$ and its representing matrix $F\in M(n\times n;%
\overline{E_{q}})$
\end{definition}

\begin{enumerate}
\item Initialize $L_{1}:=0$ and $L_{2}:=0$ and an array $v$ of length $n$ to
zero.

\item For $k$ from $1$ to $n$ do $L_{2}:=L_{2}+1$ and $v[k]:=1$ if and only
if $F_{k1}=-\infty $.

\item Compare $L_{1}$ and $L_{2}.$ If $L_{1}=L_{2}$ or $L_{2}=n$, construct
the matrix%
\begin{equation*}
F^{\prime }\in M((n-L_{2})\times (n-L_{2});E_{q})
\end{equation*}%
by deleting the $k$th row and the $k$th column of $F$ for all $k$ s.t. $%
v[k]=1.$ Then return $F^{\prime }$ and stop. If $L_{1}<L_{2}$ and $L_{2}<n$,
calculate the product $F^{\cdot 2}$, set $F:=F^{\cdot 2}$ as well as $%
L_{1}:=L_{2}$ and go to step 2.

\item If the returned matrix $F^{\prime }$ is the empty matrix ($L_{2}=n$)
we can conclude that the system $f$ is a fixed point system with $%
(0,...,0)^{t}\in \mathbf{F}_{q}^{n}$ as its unique fixed point (see Remark %
\ref{Rem.ZeroExcludedFromMonSyst.}). If $F^{\prime }$ is not the empty
matrix, the corresponding lower dimensional system $f^{\prime }:=\Psi
(F^{\prime })$ needs to be analyzed with the algorithm described in Section
4.
\end{enumerate}

Step 1 implies $n+2$ initializations. Step 2 of the algorithm requires $n$
comparisons, at most $n$ additions and at most $n$ assignments. There are $2$
comparisons in step 3. Each matrix multiplication in step 3 takes $%
2n^{3}-n^{2}$ addition or multiplication operations\footnote{%
See also the analysis of the arithmetic operations in the semiring $E_{q}$
in Section 4.} in $\overline{E_{q}}$. There is one initialization after each
matrix multiplication. The worst case scenario is given when every time the
algorithm performs step 3, the set $L_{1}$ grows by one element, forcing the
algorithm to perform $n-1$ matrix multiplications. The construction of the
matrix $F^{\prime }$ requires a number of comparisons and assignments that
is obviously bounded above by $2n^{2}.$ Summarizing, the worst case
complexity of the algorithm is bounded above by%
\begin{eqnarray*}
B(n) &:&=(n+2)+n(3n)+n2+(n-1)(2n^{3}-n^{2}+1)+2n^{2} \\
&=&2n^{4}-3n^{3}+6n^{2}+4n+1
\end{eqnarray*}%
Since%
\begin{equation*}
\lim_{n\rightarrow \infty }\frac{B(n)}{n^{4}}=2
\end{equation*}%
we can conclude $B(n)\in O($ $n^{4})$.\newline
It is pertinent to emphasize that this preprocessing algorithm represents a
primitive first attempt. Since the matrix multiplications dominate the
complexity of the algorithm, it seems meaningful to try to reduce the
complexity of the multiplication. Indeed, the rows with entries $-\infty $
are preserved during the multiplication, i.e. those rows do not need to be
calculated. In addition, if the first element of a row in the product matrix
is equal to $-\infty $ we know that all the remaining elements of that row
are going to be equal to $-\infty $ as well. As we can see, there are
possibilities of improvement. However, for the purposes of this paper, we
are satisfied with a first working algorithm of polynomial complexity.

\bibliographystyle{plain}
\bibliography{DISSMathRef}

\begin{thebibliography}{10}

\bibitem{MCA}
R.~Bartlett and M.~Garzon.
\newblock Monomial cellular automata.
\newblock {\em Complex Systems}, 7:367--388, 1993.

\bibitem{MR0299409}
Arthur~W. Burks, editor.
\newblock {\em Essays on cellular automata}.
\newblock University of Illinois Press, Urbana, Ill., 1970.

\bibitem{SDS}
H.~S.~Mortveit C.~L.~Barrett and C.~M. Reidys.
\newblock Elements of a theory of simulation ii: sequential dynamical systems.
\newblock {\em Applied Mathematics and Computation}, 107(2-3):121--136, 2000.

\bibitem{ThesisOmar}
Omar Colón-Reyes.
\newblock {\em Monomial Dyanmical Systems over Finite Fields}.
\newblock PhD thesis, Virginia Tech, Blacksburg, Virginia, 2005.

\bibitem{MR2293353}
Omar Col{\'o}n-Reyes, Abdul~Salam Jarrah, Reinhard Laubenbacher, and Bernd
  Sturmfels.
\newblock Monomial dynamical systems over finite fields.
\newblock {\em Complex Systems}, 16(4):333--342, 2006.

\bibitem{MR2112694}
Omar Col{\'o}n-Reyes, Reinhard Laubenbacher, and Bodo Pareigis.
\newblock Boolean monomial dynamical systems.
\newblock {\em Ann. Comb.}, 8(4):425--439, 2004.

\bibitem{CULL}
Paul Cull.
\newblock Linear analysis of switching nets.
\newblock {\em Kybernetik}, 8(1):31--39, 1971.

\bibitem{MR???????}
E.~Delgado-Eckert.
\newblock Canonical representatives for residue classes of a polynomial ideal
  and orthogonality.
\newblock {\em Comm. Algebra}, under review.
\newblock See temporary version at http://arxiv.org/abs/0706.1952v1.

\bibitem{Elspas}
B.~Elspas.
\newblock The theory of autonomous linear sequential networks.
\newblock {\em IRE Transactions on Circuit Theory}, CT-6:45--60, 1959.

\bibitem{MR0184874}
Frank Harary, Robert~Z. Norman, and Dorwin Cartwright.
\newblock {\em Structural models: {A}n introduction to the theory of directed
  graphs}.
\newblock John Wiley \& Sons Inc., New York, 1965.

\bibitem{MR2175374}
Ren{\'e}~A. Hern{\'a}ndez~Toledo.
\newblock Linear finite dynamical systems.
\newblock {\em Comm. Algebra}, 33(9):2977--2989, 2005.

\bibitem{JUST}
Winfried Just.
\newblock The steady state system problem is np-hard even for monotone
  quadratic boolean dynamical systems.
\newblock {\em http://www.math.ohiou.edu/~just/PAPERS/monNPh14.pdf}, 2006.

\bibitem{Kaplan}
M.~Kaplan.
\newblock {\em Computeralgebra}.
\newblock Springer-Verlag, Berlin, Heidelberg, 2005.

\bibitem{CA}
Jarkko Kari.
\newblock Theory of cellular automata: A survey.
\newblock {\em Theoretical Computer Science}, 334(1-3):3--333, 2005.

\bibitem{p-adic}
A.~Khrennikov and M.~Nilsson.
\newblock On the number of cycles of p-adic dynamical systems.
\newblock {\em Journal of Number Theory}, 90:255--264, 2001.

\bibitem{MR1429394}
R.~Lidl and H.~Niederreiter.
\newblock {\em Finite fields}, volume~20 of {\em Encyclopedia of Mathematics
  and its Applications}.
\newblock Cambridge University Press, Cambridge, second edition, 1997.
\newblock With a foreword by P. M.\ Cohn.

\bibitem{AffineSys}
D.~K. Milligan and M.~J.~D. Wilson.
\newblock The behavior of affine boolean sequential networks.
\newblock {\em Connection Sciences}, 5(2):153--167, 1993.

\bibitem{Fuzzy}
M.~Nilsson.
\newblock Fuzzy cycles in monomial dynamical systems.
\newblock {\em Far East Journal of Dynamical Systems}, 5:149--173, 2003.

\bibitem{REGER&SCHMIDT}
J.~Reger and K.~Schmidt.
\newblock A finite field framework for modeling, analysis and control of finite
  state automata.
\newblock {\em Mathematical and Computer Modelling of Dynamical Systems
  (MCMDS)}, 10(3):253--285, 2004.

\bibitem{FiniteStateAut}
J.~Reger and K.~Schmidt.
\newblock Modeling and analyzing finite state automata in the finite field f2.
\newblock {\em Mathematics and Computers in Simulation}, 66(2-3):193--206,
  2004.

\bibitem{Vasiga}
T.~Vasiga and J.~Shallit.
\newblock On the iteration of certain quadratic maps over gf(p).
\newblock {\em Discrete Mathematics}, 277:219--240, 2004.

\end{thebibliography}

\end{document}